\documentclass{csri22}

\usepackage{amsfonts,amsmath,graphicx,subfigure}
\usepackage{boldtensors}
\usepackage{bm}
\usepackage{color}
\usepackage{multirow}
\usepackage[table]{xcolor}

\newcommand{\op}[1]{\operatorname{\text{#1}}}

\newcommand{\Div}{\op{Div}}

\newcommand{\cV}{\set{V}}
\newcommand{\set}[1]{\mathcal{#1}}

\newcommand*\diff{\mathop{}\!\mathrm{d}}

\addcontentsline{toc}{chapter}{THE SCHWARZ ALTERNATING METHOD FOR THE SEAMLESS COUPLING OF NONLINEAR REDUCED ORDER MODELS AND FULL ORDER MODELS\\
{\em J. \ Barnett, I. \ Tezaur and A. \ Mota}}

\title{THE SCHWARZ ALTERNATING METHOD FOR THE SEAMLESS COUPLING OF NONLINEAR REDUCED ORDER MODELS AND FULL ORDER MODELS}

\author{Joshua Barnett\thanks{Department of Mechanical Engineering, Stanford University, jb0@stanford.edu} \and Irina Tezaur\and Alejandro Mota\thanks{Sandia National Laboratories,
ikalash@sandia.gov}}

\begin{document}

\maketitle

\begin{abstract}
Projection-based model order reduction allows for the parsimonious representation of full order models (FOMs), typically obtained through the discretization of a set of partial differential equations (PDEs) using conventional techniques (e.g., finite element, finite volume, finite difference methods) where the discretization may contain a very large number of degrees of freedom. As a consequence of this more compact representation, the resulting projection-based reduced order models (ROMs) can achieve considerable computational speedups, which are especially useful in real-time or multi-query analyses. 
	One known deficiency of projection-based ROMs is that they can suffer from a lack of robustness, stability and accuracy, especially in the predictive regime, which ultimately limits their useful application.  Another research gap that has prevented the widespread adoption of ROMs within the modeling and simulation community is the lack of theoretical and algorithmic foundations necessary for the ``plug-and-play" integration of these models into existing multi-scale and multi-physics frameworks.
		This paper describes a new methodology that has the potential to 
		address both of the aforementioned deficiencies 
		by coupling projection-based ROMs with each other as well as with conventional FOMs by means of the Schwarz alternating method \cite{JLB:schwarz1870ueber}.  
		Leveraging recent work that adapted the Schwarz alternating method to enable consistent and concurrent multi-scale coupling of finite element FOMs in solid mechanics \cite{JLB:mota2017schwarz, JLB:mota2022schwarz}, we present a new extension of the Schwarz framework that enables FOM-ROM and ROM-ROM coupling, following a domain decomposition of the physical geometry on which a PDE is posed.  In order to maintain efficiency and achieve computation speed-ups, we employ hyper-reduction via the Energy-Conserving Sampling and Weighting (ECSW) approach \cite{JLB:farhat2015structure}.  We evaluate the proposed coupling approach in the reproductive
		as well as in the predictive regime on a canonical test case that involves the dynamic propagation of a traveling wave in a nonlinear hyper-elastic material.

\end{abstract}

\section{Introduction} \label{sec:intro}

Projection-based model order reduction is a promising data-driven strategy for reducing the computational complexity 
of numerical simulations by restricting the search of the solution to a low-dimensional space spanned by a reduced basis
constructed from a limited number of high-fidelity simulations and/or physical experiments/observations.
While recent years have seen extensive investments in the development of projection-based reduced order models (ROMs)
and other data-driven models, these models are 
known to suffer from a lack of robustness, stability and accuracy, especially in the predictive regime.  Moreover,
a unified and rigorous theory for integrating these models in a ``plug-and-play" fashion
into existing multi-scale and multi-physics 
coupling frameworks (e.g., the Department of Energy's Energy Exascale Earth System Model (E3SM) \cite{JLB:E3SM}) is lacking at the present time. 

This paper presents and evaluates an approach aimed at addressing both of the aforementioned shortcomings
by advancing the Schwarz alternating method \cite{JLB:schwarz1870ueber} as a mechanism for coupling together a variety of models, including full order finite element models and projection-based ROMs constructed
using the proper orthogonal decomposition (POD)/Galerkin method.  
The Schwarz alternating method is based on the simple idea that if the solution to a partial differential equation (PDE) is known in two or more regularly-shaped domains comprising a more complex domain, these local solutions can be used to iteratively build a solution on the more complex domain.  Our coupling approach thus consists of several
ingredients, namely: (1) the decomposition of the physical domain into overlapping or non-overlapping
subdomains, which can be discretized in space by disparate meshes and in time by different time-integration schemes with different time-steps, 
(2) the definition of 
transmission boundary conditions on subdomain boundaries, and (3) the iterative solution of a sequence of subdomain problems
in which information propagates between the subdomains through the aforementioned transmission conditions.  
Without loss of generality, we develop and prototype the method in the context of a generic transient dynamic solid mechanics problem, defined by an arbitrary constitutive model embedded
within the governing PDEs.
Following an overlapping or non-overlapping
domain decomposition (DD) of the underlying geometry, we use the POD/Galerkin method with Energy-Conserving 
Sampling and Weighting (ECSW)-based hyper-reduction \cite{JLB:farhat2015structure} to reduce the problem in one or more subdomains.  
We then employ 
the Schwarz alternating method to couple the resulting subdomain ROMs with each other or with finite element-based full order models
(FOMs) in neighboring subdomains. We demonstrate that a careful formulation and implementation of the 
 transmission conditions in the ROMs being coupled is essential to the coupling method.

The methodology described in this paper is related to several existing coupling approaches developed in recent years.
First, while this work is a direct extension of the recently-developed Schwarz-based methodology for concurrent multi-scale 
FOM-FOM coupling in solid mechanics \cite{JLB:mota2017schwarz, JLB:mota2022schwarz}, it includes a number of advancements, including the extension 
of the coupling framework to: (1) FOM-ROM coupling, (2) ROM-ROM coupling, and (3) non-overlapping 
subdomains. Among the earliest authors to develop an iterative Schwarz-based 
DD approach for coupling FOMs with ROMs are Buffoni \textit{et al.} \cite{JLB:Buffoni}. 
The approach in \cite{JLB:Buffoni}
is unlike ours in that attention 
is restricted to Galerkin-free
POD ROMs, developed for the Laplace equation and the compressible Euler equations.  Other authors to consider Galerkin-free FOM-ROM and ROM-ROM
couplings are Cinquegrana \textit{et al.} \cite{JLB:Cinquegrana}
and Bergmann \textit{et al.} \cite{JLB:Bergmann2018}.
The former approach \cite{JLB:Cinquegrana} focuses on overlapping DD in the context of a 
Schwarz-like iteration scheme, but, unlike our approach, requires matching meshes at the 
subdomain interfaces.  The latter approach \cite{JLB:Bergmann2018}, termed zonal Galerkin-free POD, 
defines a minimization problem to minimize the difference between the POD reconstruction and its 
corresponding FOM solution in the overlapping region between a ROM and a FOM domain, and
is developed/investigated in the context of an unsteady flow and aerodynamic shape optimization.
While the method developed \cite{JLB:Bergmann2018} 
is not based on the Schwarz alternating formulation, the recent related work \cite{JLB:Iollo}
by Iollo \textit{et al.} demonstrates that a similar optimization-based 
coupling scheme is equivalent to an overlapping alternating Schwarz iteration
for the case of linear elliptic PDE.  
A true POD-Greedy/Galerkin non-overlapping Schwarz method for the coupling of projection-based ROMs 
developed for the specific case of symmetric elliptic PDEs is presented 
by Maier \textit{et al.} in \cite{JLB:Maier}.  

While the focus herein is restricted to projection-based ROMs, it is worth noting that the Schwarz 
alternating method has recently been extended to the case of coupling Physics-Informed
Neural Networks (PINNs) to each other following a DD in \cite{JLB:LiD3M, JLB:LiDeepDDM}.
The methods proposed in these works, termed D3M \cite{JLB:LiD3M} and DeepDDM \cite{JLB:LiDeepDDM},
inherit the benefits of DD-based ROM-ROM couplings, 
but are developed primarily for the purpose of improving the 
efficiency of the neural network training process and reducing the risk of overfitting, 
both of which are due to the global nature of the neural network ``basis functions".  

We end our literature overview by remarking that a number of non-Schwarz-based ROM-ROM and/or FOM-ROM coupling
methods have been developed in recent years, including \cite{JLB:Baiges, JLB:Wicke, JLB:Hoang,
JLB:Iapichino, JLB:LeGresley, JLB:Lucia, JLB:Maday, JLB:Corigliano1, JLB:Corigliano2,
JLB:Kerfriden1, JLB:Kerfriden2, JLB:Radermacher, JLB:deCastro2022, JLB:Smetana2022, JLB:Riffaud2021, JLB:Huang2022}.  The majority 
of these approaches are based on Lagrange multiplier or flux matching coupling formulations, and focus on either simple linear elliptic PDEs or fluid problems. We omit a detailed
assessment of these references from this paper for the sake of brevity.

The remainder of this paper is organized as follows.  In Section \ref{sec:sm_prob}, we provide the 
variational formulation of the generic solid dynamics problem considered herein, and describe its spatio-temporal
discretization.  Section \ref{sec:rom} details our nonlinear model reduction methodology for this problem, which relies on the
POD/Galerkin approach for model 
reduction and the ECSW  method \cite{JLB:farhat2015structure} for hyper-reduction.
We describe the Schwarz alternating method for 
FOM-FOM, FOM-ROM and ROM-ROM coupling in Section \ref{sec:schwarz}.  
In Section \ref{sec:results}, we evaluate the performance of the proposed Schwarz-based coupling methodology
on a problem involving dynamic wave propagation in a one-dimensional (1D) hyper-elastic bar whose material 
properties are described by the nonlinear Henky constitutive model, characterized by a logarithmic strain tensor \cite{JLB:Henky1931}.  We conclude this paper with a 
summary and a discussion of some future research directions (Section \ref{sec:conc}). 

\section{Solid mechanics problem formulation} \label{sec:sm_prob} 

Consider the Euler-Lagrange equations for a generic dynamic solid mechanics
problem in its strong form:
\begin{equation} \label{eq:euler-lagrange}
\Div ~P + \rho_0 ~B
=
\rho_0 \ddot{~\varphi}
\quad \text{in} \quad
\Omega \times I.
\end{equation}
In \eqref{eq:euler-lagrange}, $\Omega \in \mathbb{R}^d$ for $d \in \{
1,2,3\}$ is an open bounded domain, $I :=\{ t \in [0, T]\}$ is a
closed time interval with $T > 0$, and $~x = {~\varphi}(~X, t):
\Omega \times I \to \mathbb{R}^d$ is a mapping, with $~X \in \Omega$
and $t \in I$.  The symbol $~P$ denotes the first Piola-Kirchhoff stress and
$\rho_0 ~B: \Omega \to \mathbb{R}^d$ is the body force, with $\rho_0$
denoting the mass density in the reference configuration.  The
over-dot notation denotes differentiation in time, so that
$\dot{~\varphi}:= \frac{\partial ~\varphi}{\partial t}$ and
$\ddot{~\varphi}:= \frac{\partial^2 ~\varphi}{\partial t^2}$. Embedded
within $~P$ is a constitutive model, which can range from a simple
linear elastic model to a complex micro-structure model, e.g., that of
crystal plasticity.  Herein, we focus on nonlinear hyper-elastic constitutive 
models such as the Henky model \cite{JLB:Henky1931}.  The details of this 
model are provided in Section \ref{sec:results}.

Suppose that we have the following initial and boundary conditions for the
PDEs \eqref{eq:euler-lagrange}:
\begin{equation} \label{bcs}
\begin{array}{c}
  ~\varphi(~X, t_0)
   =
  ~X_0, \hspace{0.2cm}
  \dot{~\varphi}(~X, t_0)
   =
  ~v_0
   \text{  in  }
  \Omega, \\
  ~\varphi(~X, t)
   =
  ~\chi
   \text{  on  }
  \partial \Omega_{~\varphi}  \times I,
    \hspace{0.2cm} ~P ~N
   =
  ~T
   \text{  on  }
  \partial \Omega_{~T} \times I.
\end{array}
\end{equation}
In \eqref{bcs}, it is assumed the outer boundary $\partial \Omega$ is decomposed
into a Dirichlet and traction portion, $\partial \Omega_{~\varphi}$ and
$\partial \Omega_{~T}$, respectively, with $\partial \Omega = \partial
\Omega_{~\varphi} \cup \partial \Omega_{~T}$ and $\partial \Omega_{~\varphi}
\cap \partial \Omega_{~T} = \emptyset$. The prescribed boundary positions or
Dirichlet boundary conditions are $~\chi : \partial \Omega_{~\varphi} \times I
\to \mathbb{R}^3$. The symbol $~N$ denotes the unit normal on $\partial
\Omega_{~T}$.  In this work, we will assume without loss of generality that $~\chi$ is not changing in time.

It is straightforward to show that the weak variational form of
\eqref{eq:euler-lagrange} with initial and boundary conditions
\eqref{bcs} is
\begin{equation} \label{eq:var}
  \int_{I} \left[
    \int_{\Omega} \left(
      \Div ~P + \rho_0 ~B - \rho_0 \ddot{~\varphi}
    \right) \cdot ~\xi \ \diff V
    +
    \int_{\partial_{~T} \Omega} ~T \cdot ~\xi \ \diff S
  \right] \ \diff t
  =
  0,
\end{equation}
where $~\xi$ is a test function in $ \cV :=
\left\{
~\xi \in W_2^1(\Omega \times I)
:
~\xi = ~0 \text{ on } \partial_{~\varphi} \Omega \times I
\cup
\Omega \times t_0
\cup \right.$
$\left.
\Omega \times t_1
\right\}$.

Discretizing the variational form \eqref{eq:var} in space using the classical
Galerkin finite element method (FEM) \cite{JLB:HughesBook} yields the following
semi-discrete matrix problem:
\begin{equation} \label{eq:semidiscrete}
	~M \ddot{~u} + ~f^{\text{int}}(~u, \dot{~u})= ~f^{\text{ext}}.
\end{equation}
In \eqref{eq:semidiscrete}, $~M$ denotes the mass matrix, $~u := ~\varphi(~X, t)
- ~X$ is the displacement, $\ddot{~u}$ is the acceleration (also denoted by $~a$),  $~f^{\text{ext}}$ is a vector of applied external
forces, and  $~f^{\text{int}}(~u, \dot{~u})$ is the vector of internal forces due to
mechanical and other effects inside the material, where $\dot{~u}$ (also denoted by $~v$)
is the velocity.  
In the present work, 
the semi-discrete equation \eqref{eq:semidiscrete}
is advanced forward in time using the Newmark-$\beta$ time-integration scheme \cite{JLB:Newmark}.  
We will assume the FOM \eqref{eq:semidiscrete} has size $N \in \mathbb{N}$, that is, $~u \in \mathbb{R}^N$.  For convenience, 
in subsequent discussion, we transform the Dirichlet boundary 
condition \eqref{bcs} for the position into a Dirichlet boundary for the displacement, so that the boundary condition imposed on  
$\partial \Omega_{~\varphi} \times I$ is
$~u= ~u_D$, with $~u_D$ being independent of time.

\section{Model order reduction} \label{sec:rom}

Projection-based model order reduction is a 
promising, physics-based technique for reducing the computational cost associated with 
high-fidelity models such as \eqref{eq:semidiscrete}.  The basic workflow for building a projection-based ROM for a 
generic nonlinear semi-discrete problem of the form \eqref{eq:semidiscrete}
consists of three steps: (1) calculation of a reduced basis, (2) projection of the governing 
equations onto the reduced basis, and (3) hyper-reduction of the nonlinear terms in the 
projected equations.  Herein, we employ the POD \cite{JLB:Sirovich1987, JLB:Holmes1996} 
for the reduced basis generation step (step 1), the Galerkin projection method 
for the projection step (step 2), and the ECSW method \cite{JLB:farhat2015structure}
for the hyper-reduction step (step 3).  Each of these steps is described succinctly below.

\subsection{Proper orthogonal decomposition (POD)} \label{sec:pod}

The POD is a mathematical procedure that, given an ensemble of
data and an inner product, constructs a basis for the ensemble that is
optimal in the sense that it describes more energy (on average) of the ensemble in the chosen inner
product than any other linear basis of the same dimension $M$. 
The ensemble ${~w^s \in \mathbb{R}^N: s = 1, . . . , S}$ 
is typically a set of $S$ instantaneous snapshots of a numerical solution field, 
collected for $S$ values
of a parameter of interest, and/or at $S$ different times. 
For solid mechanics problems, a natural choice 
for the snapshots is $~w^s = ~u^s$, where the ensemble $\{~u^s\}$ denotes a set of snapshots for the displacement field.  It is noted that one can use in place of on addition to $\{~u^s\}$
 snapshots of the velocity ($\{~v^s\}$) and/or acceleration ($\{~a^s\}$) fields.

Following the so-called 
``method of snapshots" \cite{JLB:Sirovich1987, JLB:Holmes1996}, a POD basis $~\Phi_M \in \mathbb{R}^{N \times M}$ 
of dimension $M$ is obtained by performing a singular 
value decomposition (SVD) of a snapshot matrix $~W := \left[ ~w^1, ..., ~w^S\right] \in \mathbb{R}^{N\times S}$
such that $~W = ~\Phi ~\Sigma ~V^T$ and defining $~\Phi_M$ as the matrix containing the first $M$ columns of $~\Phi$. 
Letting
\begin{equation} \label{eq:energy_pod}
	\mathcal{E}_{\text{POD}}(M) := \frac{\sum_{i=1}^M \sigma_i^2}{\sum_{i=1}^S \sigma_i^2}
\end{equation}
denote the energy associated with a POD basis of size $M$, where $\sigma_i$ denotes
the $i^{th}$ singular value of $~W$, the basis size $M$ is typically 
selected based on an energy criterion, where $100\mathcal{E}_{\text{POD}}(M)$ is the percent
energy captured by a given POD basis. 

As discussed in Section \ref{sec:proj}, it is often desirable to construct the POD basis $~\Phi_M$
such that this basis satisfies homogeneous Dirichlet boundary conditions at some 
pre-defined indices $~i_d \in \mathbb{N}^d$ where $ d < N$.  It is straightforward 
to accomplish this by simply zeroing out the snapshots at the Dirichlet degrees of freedom (dofs)
prior to calculating the SVD, that is, setting $~w^i(~i_d) = ~0$ for $i=1, ..., S$.


\subsection{Galerkin projection} \label{sec:proj}
As discussed in \cite{JLB:farhat2015structure, JLB:Tiso2013}, Galerkin projection 
is in general the method of choice for solid mechanics and structural dynamics problems, as it preserves the Hamiltonian structure of the underlying system of PDEs \cite{JLB:Lall_PhysicaD_2003}.  
The method starts by approximating the FOM displacement solution to \eqref{eq:semidiscrete} as
\begin{equation} \label{eq:disp_decomp}
	~u \approx \tilde{~u} = \bar{~u} + ~\Phi_M \hat{~u},
\end{equation}
where $\hat{~u} \in \mathbb{R}^M$ is the vector of unknown modal amplitudes, to be solved for in the ROM and
$\bar{~u}\in \mathbb{R}^N$ is a (possibly time-dependent) reference state.  Similar approximations
can be made for the velocity and acceleration fields, $~v:=\dot{~u}$ 
and $~a := \ddot{~u}$, respectively, namely:
\begin{equation} \label{eq:acce_decomp}
	\begin{array}{c}
	~v \approx \tilde{~v} = \bar{~v} + ~\Phi_M \hat{~v},\\
	~a \approx \tilde{~a} = \bar{~a} + ~\Phi_M \hat{~a}.
	\end{array}
\end{equation}
Here, $\bar{~v}$ and $\bar{~a}$ are defined analogously to $\bar{~u}$, 
and similarly for $\hat{~v}$ and $\hat{~a}$.  

In the present formulation,
the reference states $\bar{~u}$, $\bar{~v}$ and $\bar{~a}$ 
are used to prescribe strongly Dirichlet boundary conditions 
within the ROM for the displacement, velocity and acceleration fields, respectively.  
This is done following the approach of Gunzburger \textit{et al.} \cite{JLB:Gunzburger2007}.  
Suppose the Dirichlet 
boundary conditions of interest are $~u(~i_d) = ~u_D$, $~v(~i_d) = ~v_D$ 
and $~a(~i_d) = ~a_D$.  Let $~i_u$ denote the unconstrained indices
at which the solution to \eqref{eq:semidiscrete} is sought.
It is straightforward to 
see from \eqref{eq:disp_decomp}--\eqref{eq:acce_decomp} that, 
if the POD modes $~\Phi_M = \left[~\phi_1, ..., ~\phi_M\right]$ are calculated such that
$~\phi_i(~i_d) = ~0$ for $i = 1, ..., M$ and 
\begin{equation}
\begin{array}{ll}
	\bar{~u}(~i_d) = ~u_D, & \bar{~u}(~i_u) = ~0, \\
	\bar{~v}(~i_d) = ~v_D, & \bar{~v}(~i_u) = ~0, \\
	\bar{~a}(~i_d) = ~a_D, & \bar{~a}(~i_u) = ~0, \\
\end{array}
\end{equation}
the ROM solution will satisfy the prescribed Dirichlet boundary conditions
in the strong sense.  

The ROM for \eqref{eq:semidiscrete} is obtained by substituting the decompositions \eqref{eq:disp_decomp}--\eqref{eq:acce_decomp}
into the FOM equations \eqref{eq:semidiscrete}, and projecting these equations onto the reduced basis $~\Phi_M$.  It is straightforward to verify that 
doing this and moving all Dirichlet dofs to the right-hand side of the resulting 
system of equations yields a semi-discrete problem of the form
\begin{equation} \label{eq:rom}
	\hat{~M}_{uu}\hat{~a} + \hat{~f}^{\text{int}}_u(\tilde{~u}, \tilde{~v}) = \hat{~f}^{\text{ext}}_u,
\end{equation}
where 
\begin{equation} \label{eq:defns}
\begin{array}{c}
	\hat{~M}_{uu} := ~\Phi_{M,u}^T ~M_{uu} ~\Phi_{M,u}, \hspace{0.1cm} 
	~f^{\text{int}}_u(\tilde{~u}, \tilde{~v}):=~\Phi_{M,u}^T ~f_u^{\text{int}}(\tilde{~u}, \tilde{~v}),\\ 
	~f^{\text{ext}}_u := ~\Phi_{M,u}^T (~f_u^{\text{ext}} - ~M_{ud} \bar{~a}_d).
\end{array}
\end{equation}
In \eqref{eq:defns}, $~\Phi_{M,u}:=~\Phi_M(~i_u, :)$, $~M_{uu}:=~M(~i_u, ~i_u)$, $~f_u^{\text{int}}:= ~f^{\text{int}}(~i_u)$, 
$~f^{\text{ext}}_u:=~f^{\text{ext}}(~i_u)$, $~M_{ud}:=~M(~i_u, ~i_d)$
and $\bar{~a}_d:=\bar{~a}(~i_d)$.   \\

\noindent {\bf \textit{Remark 1.}} It is noted that, while all three variables $\tilde{~u}$, $\tilde{~v}$ and $\tilde{~a}$ appear in the ROM system \eqref{eq:defns}, we are not considering the displacement, velocity and acceleration variables as independent fields with separate generalized coordinates in our ROM construction approach.  In particular, $\hat{~v}:=\frac{d\hat{~u}}{dt}$ and 
$\hat{~a}:=\frac{d^2\hat{~u}}{dt^2}$ in \eqref{eq:acce_decomp}.  Taking this approach ensures consistency between the displacement, velocity and acceleration solutions computed within the ROM.

\subsection{Hyper-reduction via ECSW} \label{sec:ecsw}

The POD/Galerkin approach to model reduction described in Sections \ref{sec:pod} and \ref{sec:proj}
is not efficient for nonlinear problems, as the projection of the nonlinear
terms appearing in the ROM system \eqref{eq:rom}
requires algebraic operations that scale with the dimension of the original full order model $N$.  This
problem can be circumvented through the use of a procedure known as hyper-reduction.  Here, we 
rely on a specific ``project-then-approximate" hyper-reduction approach known as ECSW \cite{JLB:farhat2015structure}.
We select this approach, as it has been shown in \cite{JLB:farhat2015structure}
to preserve the Lagrangian structure associated 
with Hamilton's principle for second-order dynamical systems of the form \eqref{eq:semidiscrete}.
The resulting hyper-reduced ROM, termed HROM, is thus able to preserve the numerical stability properties of the 
Newmark-$\beta$ time-integration scheme applied to advance the ROM system \eqref{eq:rom} forward in time.

In the context of solid mechanics, ECSW can be described as cubature-based approach with the goal of accurately
estimating the projected internal forcing term (the ``project" part of ``project-then-approximate") 
$\hat{~f}^{\text{int}}_u(\tilde{~u}, \tilde{~v})$
in \eqref{eq:rom} as opposed to
directly estimating the nonlinear forcing term and then projecting the result with the appropriate basis. This approximation can be rewritten as a summation of each of the element-wise contributions

\begin{equation}
    \label{eq:ecsw_element-wise}
    \begin{array}{rl}
        \hat{~f}^{\text{int}}_u(\tilde{~u}, \tilde{~v}) & = \sum_{e \in \mathcal{E}} ~\Phi_{M}^T ~L_e^T ~f{_{u, e}^{\text{int}}}(~L_e\tilde{~u}, ~L_e\tilde{~v}),\\
        \hat{~K}(\tilde{~u}, \tilde{~v}) & =  \sum_{e \in \mathcal{E}} ~\Phi_{M}^T ~L_e^T ~K_e(~L_e\tilde{~u}, ~L_e\tilde{~v}) ~L_e ~\Phi_{M}.
    \end{array}
\end{equation}
In \eqref{eq:ecsw_element-wise}:  
\begin{itemize}
    \item $~f{_{u, e}^{\text{int}}}$ is the contribution to the internal forcing vector due to mesh element $e$;
    \item $~K$, $\hat{~K}$ and $~K_e$ are the tangent stiffness matrix, reduced tangent stiffness matrix, and contribution to the tangent stiffness matrix due to element $e$ respectively, which are formed for use in implicit time-stepping methods;
    \item $\mathcal{E} = \{e_1, e_2, \dots, e_N\}$ is the set of all mesh elements $e_i$ given by the finite element discretization of the FOM; and 
    \item $~L_e \in \{0, 1\}^{n_d^e \times N}$ is a Boolean matrix selecting the $n_d^e$ degrees of freedom associated with mesh element $e$.
\end{itemize}

The next step is to relax the equality in \eqref{eq:ecsw_element-wise}  such that we only approximate the reduced internal forcing vector (or tangent stiffness matrix) by sampling a subset of the mesh elements and weighting them appropriately. This can be written as 

\begin{equation}
    \label{eq:ecsw_element-wise_approx}
    \begin{array}{rl}
        \hat{~f}^{\text{int}}_u(\tilde{~u}, \tilde{~v}) & \approx \sum_{e \in \tilde{\mathcal{E}}} \xi_e ~\Phi_{M}^T ~L_e^T ~f{_{u, e}^{\text{int}}}(~L_e\tilde{~u}, ~L_e\tilde{~v}),\\
        \hat{~K}(\tilde{~u}, \tilde{~v}) & \approx  \sum_{e \in \tilde{\mathcal{E}}} \xi_e ~\Phi_{M}^T ~L_e^T ~K_e(~L_e\tilde{~u}, ~L_e\tilde{~v}) ~L_e ~\Phi_{M},
    \end{array}
\end{equation}


\noindent where  $\tilde{\mathcal{E}} \subset \mathcal{E}$ such that $|\tilde{\mathcal{E}}| = N_e \leq N$ and $\xi_e \in \mathbb{R}_{\geq 0}$ is the positive weight associated with sampled mesh element $e \in \tilde{\mathcal{E}}$. In order to determine both the sampled mesh elements and their weights, a system of matrices and vectors is formed using $N_h$ snapshots of $~u^s$ and $~v^s$ -- the same snapshots used to compute the POD modes $~\Phi_M$.  Using these snapshots, we form the following quantities:

\begin{equation}
    \label{eq:ecsw_solve}
    \begin{array}{rl}
        \tilde{~u}^s & = ~\Phi_{M} ~\Phi_{M}^T \left(~u^s - \bar{~u} \right) + \bar{~u}, \\
        \tilde{~v}^s & = ~\Phi_{M} ~\Phi_{M}^T \left(~v^s - \bar{~v} \right) + \bar{~v}, \\
        c_{se} & = ~\Phi_{M}^T ~L_e^T ~f{_{u, e}^{\text{int}}}(~L_e \tilde{~u}^s, ~L_e \tilde{~v}^s), \\
        d_s & = ~\Phi_{M}^T ~f{_{u, e}^{\text{int}}}(\tilde{~u}^s, \tilde{~v}^s),
    \end{array}
\end{equation}
where  
$\tilde{~u}^s$ and $\tilde{~v}^s$ are the reconstruction of the FOM snapshots using the affine approximations in \eqref{eq:disp_decomp} and \eqref{eq:acce_decomp}.

Given the quantities defined above, we can 
construct a system of equations 
\begin{equation} \label{eq:ecsw_lin_system}
    ~C~\xi = ~d, \hspace{0.2cm} ~\xi \in \mathbb{R}^N_{\geq 0},
\end{equation}
where
\begin{equation*}
    ~C
    := \begin{pmatrix}
    c_{11} & \dots & c_{1 N} \\
    \vdots & \ddots & \vdots \\
    c_{N_h 1} & \dots & c_{N_h N}
    \end{pmatrix} \in \mathbb{R}^{M N_h \times N}, \quad
    ~d := \begin{pmatrix}
        d_1 \\ \vdots \\ d_{N_h}
    \end{pmatrix} \in \mathbb{R}^{M N_h}
\end{equation*}
Equation \eqref{eq:ecsw_lin_system} 
is satisfied by the choice $~\xi = ~1$. This choice results in a reduced set of element equivalent to the FOM, providing no increase in computational efficiency because $N_e$, defined as the number of non-zero weights in $~\xi$ sampling the finite element mesh, is the same as the number of mesh elements. Instead, the equality in the equation \eqref{eq:ecsw_lin_system} is relaxed to instead solve a constrained optimization problem
\begin{equation}
    \label{eq:ecsw_constrained_opt}
    ~\xi = \arg \min_{~x \in \mathbb{R}^N} \left|\left| ~C ~x - ~d\right|\right|_2 \text{ subject to } ~x \geq ~0,
\end{equation}
which, for the purposes of the present work, is solved using MATLAB's \texttt{lsqnonneg} function, with an early termination criterion with a solution step size tolerance of $10^{-4}$. By relaxing the equality and including an early termination condition, the system can be solved approximately with a sparse solution for $~\xi$ such that, ideally, $N_e \ll N$. Note that our early termination criterion differs from the typical early termination criterion, 
in which the algorithm terminates once  $||~C ~\xi - ~d ||_2 / ||~d||_2$ is less than some chosen tolerance ($<1$). 
As a result, all of the HROM results presented later in the present work have a consistent termination criterion with respect to its MATLAB implementation; however, the relative error tolerance of the selected reduced elements will differ. 
\section{The Schwarz alternating method for multi-scale coupling} \label{sec:schwarz}

\begin{figure}[htbp!]
        \begin{center}
                \subfigure[Overlapping]{
      \includegraphics[width=0.4\textwidth]
                      {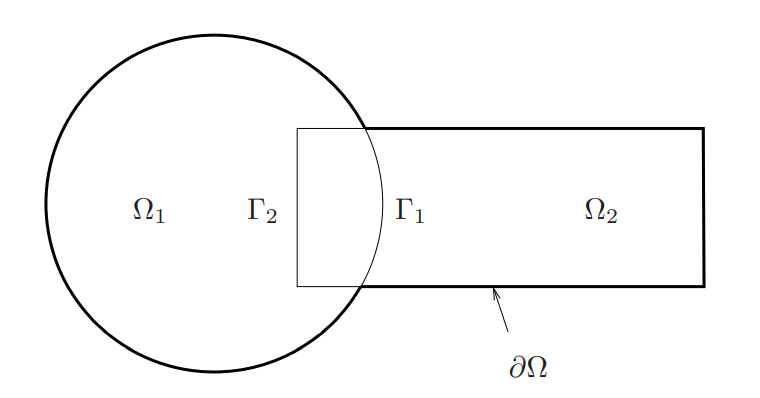}}
    \subfigure[Non-overlapping]{
      \includegraphics[width=0.4\textwidth]
                      {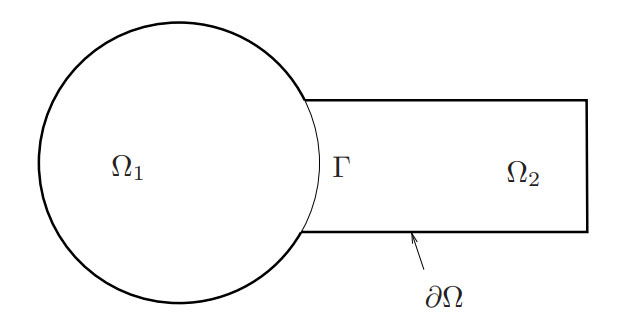}}
        \end{center}
        \caption{Illustration showing overlapping and non-overlapping domain
    decomposition.}
        \label{fig:dd}
\end{figure}

In this section, we describe two variants of the Schwarz alternating method for concurrent multiscale 
coupling.  
The first is based on an overlapping DD (Figure \ref{fig:dd}(a)) and the second
is based on a non-overlapping DD (Figure \ref{fig:dd}(b)).  Consider without loss of generality a partition of the domain 
$\Omega$ into two open subdomains $\Omega_1$ and $\Omega_2$, as shown in Figure \ref{fig:dd}, and suppose 
we are interested in applying the method to a single-physics semi-discretized PDE of the form \eqref{eq:semidiscrete}.  

\begin{figure}[htbp!]
        \begin{center}
      \includegraphics[width=1.0\textwidth]
                      {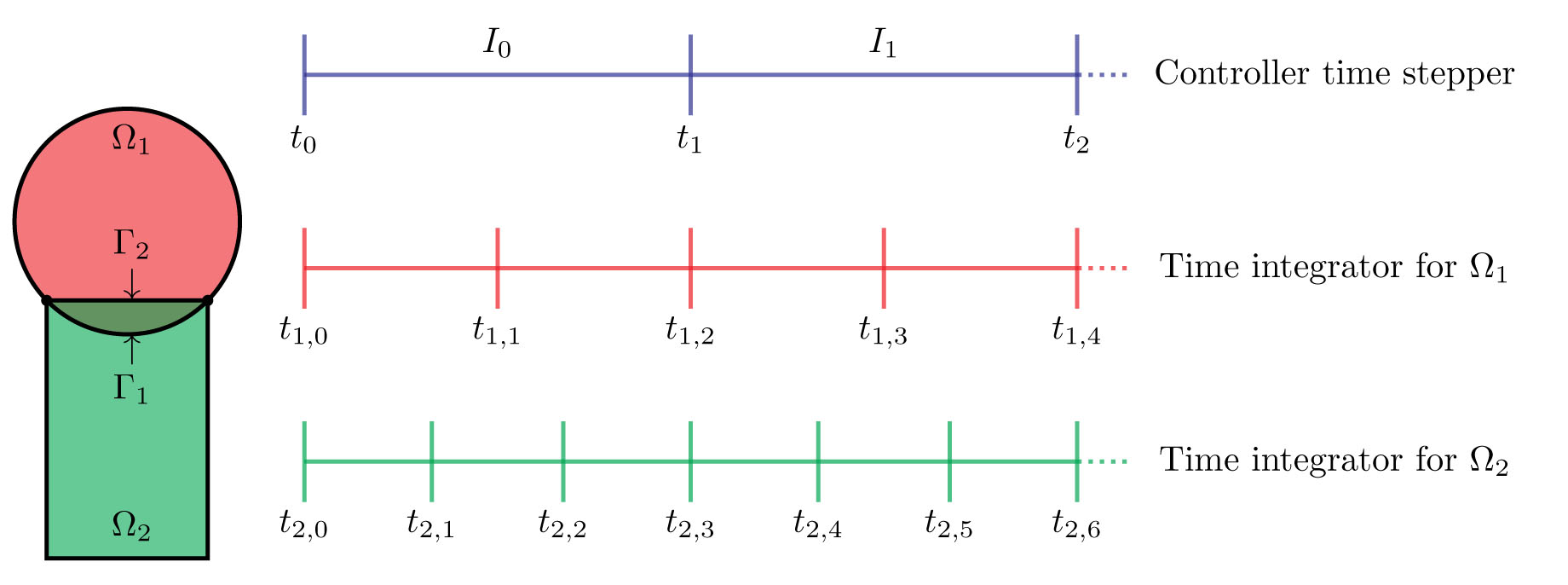}
        \end{center}
        \caption{Illustration of controller time-stepper-based time advancement of a coupled 
	problem involving two overlapping subdomains $\Omega_1$ and $\Omega_2$.}
        \label{fig:controller_time_stepper}
\end{figure}

It is important to note that, to avoid a space-time discretization for the application of the Schwarz method, it proves convenient to subdivide the time domain into ``controller time-intervals", $I_0:=[t_0, t_1]$, $I_1:=[t_1,t_2]$,..., as shown in Figure \ref{fig:controller_time_stepper}.
The controller time-intervals are convenient markers for events of interest, and for synchronization
of the Schwarz algorithm applied only in the space domain. They also define the periods or intervals in which the solutions of the initial boundary value problem \eqref{eq:semidiscrete} are determined by means of the Schwarz alternating method; effectively, the Schwarz iteration process is converged within a controller time interval $I_N$ before 
advancing to the next controller time interval $I_{N+1}$. The interested reader is referred to \cite{JLB:mota2022schwarz} and \cite{JLB:Hoy2021} for details.

\subsection{Overlapping Schwarz formulation} \label{sec:overlap}

Suppose a DD has been performed such that $\Omega_1$ and $\Omega_2$ are overlapping, that is, such that 
$\Omega_1 \cap \Omega_2 \neq \emptyset$, as shown in Figure \ref{fig:dd}(a).  In this case, the Schwarz alternating iteration
for the specific case of \eqref{eq:semidiscrete} with Dirichlet boundary conditions on $\partial \Omega$
takes the following form within each time-interval $I_N$ (Figure \ref{fig:controller_time_stepper}):
\begin{small}
\begin{equation} \label{mechanics_rr1}
        \left\{ \begin{aligned}
		~M_1 \ddot{~u}_1^{(n+1)} + ~f_1^{\text{int;(n+1)}} & = ~f_1^{\text{ext;(n+1)}}, \text{ in } \Omega_1,
      \\
                ~u_1^{(n+1)} & = ~u_D,  \text{  on  }       \partial_{~\varphi} \Omega_1 \backslash
      \Gamma_1,
      \\
            ~u_1^{(n+1)} & = ~u_2^{(n)},  \text{  on  }  \Gamma_1,\\
            \dot{~u}_1^{(n+1)} &= \dot{~u}_2^{(n)}, \text{  on  }  \Gamma_1,\\
             \ddot{~u}_1^{(n+1)} &= \ddot{~u}_2^{(n)}, \text{  on  }  \Gamma_1,
      \\
        \end{aligned} \right. \hspace{0.2cm}
        \left\{ \begin{aligned}
		~M_2 \ddot{~u}_2^{(n+1)} + ~f_2^{\text{int;(n+1)}} & = ~f_2^{\text{ext;(n+1)}}, \text{ in } \Omega_2,
      \\
      ~u_2^{(n+1)} & = ~u_D,  \text{  on  } 
      \partial_{~\varphi} \Omega_2 \backslash
      \Gamma_2,
      \\
            ~u_2^{(n+1)} & = ~u_1^{(n+1)},  \text{  on  }  \Gamma_2, \\
             \dot{~u}_2^{(n+1)} & = \dot{~u}_1^{(n+1)},  \text{  on  }  \Gamma_2, \\
              \ddot{~u}_2^{(n+1)} & = \ddot{~u}_1^{(n+1)},  \text{  on  }  \Gamma_2.
      \\
    \end{aligned} \right.
\end{equation}
\end{small}
In \eqref{mechanics_rr1}, $(n)$ for $n = 0, 1, 2, .... $ denotes the Schwarz iterations number,
and we have introduced the short-hand $~f_i^{\text{int}}:=~f_i^{\text{int}}(~u, \dot{~u})$.  The reader 
can observe that the transmission boundary conditions on the Schwarz boundaries $\Gamma_1$ and $\Gamma_2$ 
are of the Dirichlet type. 
The Schwarz Dirichlet boundary 
conditions are applied for not just the displacement field $~u$, but also for the velocity and acceleration
fields, $~v$ and $~a$, as we found in our earlier work \cite{JLB:mota2017schwarz, JLB:mota2022schwarz} that these additional conditions
are essential for maintaining accuracy within the acceleration and velocity fields when 
performing Schwarz-based coupling.
The iteration \eqref{mechanics_rr1} 
continues until convergence is reached.  Herein, convergence of the method is declared
when $||~u^{(n+1)} - ~u^{(n)}||_2/||~u^{(n)}|| < \delta$, $||~v^{(n+1)} - ~v^{(n)}||_2/||~v^{(n)}||_2 < \delta$ and $||~a^{(n+1)}-~a^{(n)}||_2/||~a^{(n)}||_2<\delta$ (where $~v^{(n)} := \dot{~u}^{(n)}$ and $~a^{(n)}:=\ddot{~u}^{(n)}$), for some specified Schwarz tolerance $\delta$.
It can be shown \cite{JLB:mota2022schwarz, JLB:Lions1988} that
the sequence defined by \eqref{mechanics_rr1} will converge to the solution of the underlying
single-domain problem provided that problem is well-posed, and the overlap is non-empty ($\Omega_1 \cap \Omega_2 \neq \emptyset$).  

\subsection{Non-overlapping Schwarz formulation} \label{sec:nonoverlap}

In certain applications, one is interested in performing coupling in a non-overlapping fashion (Figure \ref{fig:dd}(b)), e.g., 
for multi-physics applications such as fluid-structure interaction (FSI).  As noted in Section \ref{sec:overlap}, 
the overlapping formulation \eqref{mechanics_rr1} will not converge in the case the overlap 
region is empty.  This situation can be remedied by changing the transmission conditions applied 
on the Schwarz boundary $\Gamma$ shown in Figure \ref{fig:dd}(b).  Specifically, it can be shown \cite{JLB:Zanolli1987, JLB:Funaro1988} that applying the 
following alternating 
Dirichlet-Neumann (or, more-specifically, displacement-traction) iteration in each controller time-interval $I_N$ (Figure \ref{fig:controller_time_stepper}) gives rise to a convergent sequence of iterations:
\begin{small}
\begin{equation} \label{mechanics_rr2}
\begin{array}{l}
        \left\{ \begin{aligned}
		~M_1 \ddot{~u}_1^{(n+1)} + ~f_1^{\text{int;(n+1)}} & = ~f_1^{\text{ext;(n+1)}}, \text{ in }\Omega_1,
      \\
                ~u_1^{(n+1)} & = ~u_D,  \text{  on  }       \partial_{~\varphi} \Omega_1 \backslash
      \Gamma,
      \\
		~u_1^{(n+1)} & = ~\lambda_{n+1}(~u),  \text{  on  }  \Gamma, \\
			\dot{~u}_1^{(n+1)} & = ~\lambda_{n+1}(~v),  \text{  on  }  \Gamma, \\
				\ddot{~u}_1^{(n+1)} & = ~\lambda_{n+1}(~a),  \text{  on  }  \Gamma,
      \\
        \end{aligned} \right.  \hspace{0.2cm}
        \left\{ \begin{aligned}
		~M_2 \ddot{~u}_2^{(n+1)} + ~f_2^{\text{int;(n+1)}} & = ~f_2^{\text{ext;(n+1)}}, \text{ in } \Omega_2,
      \\
      ~u_2^{(n+1)} & = ~u_D,  \text{  on  }
      \partial_{~\varphi} \Omega_2 \backslash
      \Gamma,
      \\
            ~T_2^{(n+1)} & = ~T_1^{(n+1)},  \text{  on  }  \Gamma, 
      \\
    \end{aligned} \right.
    \end{array}
\end{equation}
\end{small}
where 
\begin{equation} \label{eq:relaxation}
	~\lambda_{n+1}(~u) = \theta  ~u_2^{(n)} + (1-\theta) ~\lambda_n(~u) \text{  on }\Gamma, \text{  for } n \geq 1,
\end{equation}
and similarly for $~\lambda_{n+1}(~v)$ and $~\lambda_{n+1}(~a)$, with $\theta \in (0,1]$ denoting a so-called relaxation parameter.  In our present implementation, $~\lambda_n(\cdot)$ in \eqref{eq:relaxation} is initialized to zero, i.e., $~\lambda_0(~\cdot)=~0$.    
As for the overlapping case, the iteration \eqref{mechanics_rr2} continues until convergence is reached.
If $\theta=1$ in \eqref{eq:relaxation}, 
there is no relaxation and the Schwarz boundary condition for $\Omega_1$ in \eqref{mechanics_rr2}
reduces to a regular Dirichlet boundary
condition.  
It has been demonstrated in the literature 
\cite{JLB:Zanolli1987, JLB:Funaro1988, JLB:Cote2005, JLB:Kwok2014}
that selecting $\theta < 1$ (i.e., including relaxation)
can reduce the number of Schwarz iterations required for convergence in certain applications.

It is noted that the formulation \eqref{mechanics_rr2} is not unique; in particular, it is possible to 
create a convergent Schwarz iteration for the non-overlapping DD case if one uses Robin-Robin 
transmission conditions on the Schwarz boundary $\Gamma$, as shown in \cite{JLB:LionsNonOverlap1988}.  Here, we chose 
to develop the alternating Dirichlet-Neumann formulation over the Robin-Robin formulation because
Dirichlet (displacement) and Neumann (traction) boundary conditions are readily available in most solid mechanics codes, 
unlike Robin boundary conditions. \\

\noindent {\bf \textit{Remark 2.}} Although this is not explicitly stated herein in order to avoid 
introducing overly complex notation, we note that the Schwarz boundary conditions appearing in \eqref{mechanics_rr1} 
and \eqref{mechanics_rr2} require the definition of projection operators $P_{\Omega_j \to \Gamma_i}[\cdot]$
and $P_{\Omega_j \to \Gamma}[\cdot]$, which take the solution in $\Omega_j$, and project and interpolate it onto $\Gamma_i$ or $\Gamma$.  
For more details on how this projection operator can be constructed, the interested reader is referred to
\cite{JLB:mota2017schwarz, JLB:mota2022schwarz}.  
For 1D problems, such as those considered herein, the required projection and interpolation is trivial.
It is noted that, in multiple 
spatial dimensions, the non-overlapping Schwarz formulation \eqref{mechanics_rr2} 
requires the development of operators for consistent transfer (interpolation) of traction boundary conditions using 
the concept of prolongation/restriction.  We plan to investigate
the usage of the Compadre toolkit \cite{JLB:Compadre} for this task.  \\

\noindent {\bf \textit{Remark 3.}} The Schwarz iterations \eqref{mechanics_rr1} and \eqref{mechanics_rr2}
as written are performing something that is often referred to as the multiplicative Schwarz algorithm \cite{JLB:Gander2008}, meaning 
that the Schwarz iteration in subdomain $\Omega_2$ at time-step $n+1$ depends on the Schwarz solution
in subdomain $\Omega_1$ at time-step $n+1$.  The Schwarz iteration can be modified to achieve
what is commonly referred to as the additive Schwarz method \cite{JLB:Gander2008} by applying boundary conditions from the $n^{th}$
Schwarz iteration in $\Omega_1$ to the $(n+1)^{st}$ $\Omega_2$ sub-problem.  If this change is made, the Schwarz iteration
sequences \eqref{mechanics_rr1} and \eqref{mechanics_rr2} can be parallelized over the number of subdomains.  While we do not consider the additive 
variant of the Schwarz method in the present work, preliminary results have suggested that the method
does not reduce solution accuracy and can achieve speed-ups if parallelized appropriately.

\subsection{Extension to ROM-ROM and FOM-ROM coupling} \label{sec:schwarz_rom} 

As mentioned earlier, this paper presents two novel extensions of Schwarz-based concurrent
coupling: (1) the extension of the original overlapping multi-scale coupling
formulation \cite{JLB:mota2022schwarz} to the non-overlapping DD case (Section \ref{sec:nonoverlap}), and (2)
the extension of the Schwarz coupling framework to the case when projection-based ROMs
are being coupled to each other as well as to FOMs.  

For the latter extension, it is particularly
important to be able to prescribe time-varying Dirichlet boundary conditions strongly
within a given subdomain ROM.  We achieve this by applying the approach detailed in Section \ref{sec:proj}.  
In particular, the reference states $\bar{~u}$, $\bar{~v}$ and $\bar{~a}$ in
\eqref{eq:disp_decomp}--\eqref{eq:acce_decomp} are updated in each Schwarz iteration requiring
application of a Dirichlet boundary condition on a Schwarz boundary.  Given this approach,
it is straightforward to extend the formulations \eqref{mechanics_rr1} and \eqref{mechanics_rr2} 
to projection-based ROMs, by simply replacing the FOM semi-discrete equations in these iterations
with their ROM analogs \eqref{eq:rom}. 

Equally important to our Schwarz-based coupling strategy is that the boundary nodes on which the Schwarz transmission boundary conditions are  imposed be included in the sample mesh when one is using hyper-reduction (Section \ref{sec:ecsw}). In our ROM boundary condition formulation (see Section \ref{sec:proj}), this is done automatically, as the boundary nodes are effectively removed from the ROM solve and incorporated them into the reference states $\bar{~u}$, $\bar{~v}$ and $\bar{~a}$ appearing in \eqref{eq:rom}. 

Another important detail in performing Schwarz-based coupling involving projection-based ROMs 
relates to the snapshot collection strategy.  Ideally, one would collect snapshots 
by performing simulations on the subdomains being coupled independently, that is, without running 
a coupled problem on the two (or more) domains.  Snapshot collection strategies of this sort 
(e.g., using approaches such as oversampling \cite{JLB:Smetana2022}) will be examined in future work.  In the present work, 
snapshots are generated in each subdomain by running a FOM-FOM coupled problem
using the same DD as the targeted ROM-ROM or FOM-ROM coupling, saving the converged
solutions within each subdomain, and utilizing these solutions for the construction of subdomain-local POD 
bases.  We note that one could alternatively perform a single-domain FOM simulation on the full domain $\Omega$, and use snapshots from that simulation restricted and/or interpolated onto each subdomain to generate subdomain POD bases; this second approach is not considered herein.



\section{Numerical results: nonlinear wave propagation problem} \label{sec:results}


The alternating Schwarz-based coupling strategy described in Section \ref{sec:schwarz} is evaluated on 
a 1D wave propagation problem for the solid dynamics PDEs given in Section \ref{sec:sm_prob}. 
Consider a simple beam geometry of length $1$, so that 
$\Omega = (0,1) \in \mathbb{R}$. We will assume that this geometry is clamped at both ends, 
that is, $u(0,t) = u(1,t) = 0$ for all $t \geq 0$.  The problem is initiated by prescribing 
an initial condition $u(x,0) = f(x)$ for $x \in \Omega$, where $f \in C^0(\Omega)$.  In the numerical results 
presented herein, we consider two initial conditions: 
\begin{equation} \label{eq:ic1}
	f(x) = \frac{a}{2}\exp\left[-\frac{(x-b)^2}{2s^2} \right],
\end{equation}
and
\begin{equation} \label{eq:ic2}
	f(x) = a \left[\tanh(-b(x-s)) + \tanh(b(x-1+s))\right],
\end{equation}
for $a, b, s \in \mathbb{R}$. 
We will refer to \eqref{eq:ic1} as the ``Symmetric Gaussian" initial condition,
and to \eqref{eq:ic2} as the ``Rounded Square" initial condition.
For both initial conditions, the problem is run 
until a final  time 
of $T = 1.0 \times 10^{-3}$.    

As demonstrated in Section 3.3 of \cite{JLB:mota2022schwarz}, if
$\Omega$ is assumed to be made up of a linear elastic material,
it is possible to derive an exact analytical solution to the governing
PDEs in terms of the initial condition $f(x)$ using the method of
superposition.  Since our objective herein is to evaluate the proposed
model reduction coupling methodology on a nonlinear problem, we will
assume $\Omega$ is made up of a nonlinear hyper-elastic material whose
behavior is described by the Henky constitutive model
\cite{JLB:Henky1931} with Young's modulus $E = 1$GPa and density
$\rho= 1000$ kg/m$^3$. The defining feature of Henky-type material
models is that they possess a quadratic energy density function and
that they use a logarithmic strain. For the 1D case,
this reduces to
\begin{equation} \label{Henky}
  \begin{split}
    \begin{aligned}
      \lambda(x, t) & := \frac{\partial \left[u(x, t) + x\right]}{\partial x},
      \\
      \varepsilon(x, t) & := \log \lambda(x, t),
      \\
      W(x, t) & := \frac{E}{2} \varepsilon(x, t)^2,
    \end{aligned}
  \end{split}
\end{equation}
where $\lambda$ is the stretch, $\varepsilon$ is the strain, and $W$
is the energy density. Henky-type material models are popular because
they are seen as the natural extension of linear elasticity to the
finite deformation regime, given that
\begin{equation} \label{linear-elastic}
  \begin{split}
    \begin{aligned}
      \varepsilon(x, t) & := \lambda(x, t) - 1,
      \\
      W(x, t) & := \frac{E}{2} \varepsilon(x, t)^2,
    \end{aligned}
  \end{split}
\end{equation}
for a linear elastic model, where its strain may be interpreted as an
approximation to the logarithmic strain of the Henky model. 
Since the Henky material model
variant of the problem of interest does not have an exact analytical solution, 
we began our study by first verifying that convergence of the single-domain FOM
solution to the problem with respect to both mesh and time-step refinement is
achieved.  
The details of this study are omitted herein for the sake of brevity.

Several snapshots of the solution to the Symmetric Gaussian version of our model 
problem are plotted in Figure \ref{fig:repro_solns}.  
The reader can observe that the initial condition splits in two and propagates 
to both the left and the right, eventually reflecting back from the clamped boundaries.
The behavior of the solution to the Rounded Square version of our model problem is similar.
For both problem variants, a sharp gradient forms within the acceleration field
over time.  Although simple and 1D, the nonlinear 
wave propagation problem is a challenging test case for coupling methods, as artifacts 
introduced by these schemes will be clearly evident in the numerical solution.
We additionally remark that our model problem is particularly difficult for ROMs, which 
are known to struggle on traveling wave problems and problems with sharp gradients.
Here, we explore the possibility of mitigating
this difficulty by: (1) restricting the ROM to a part of the domain and coupling it to a FOM in the remainder of the domain, as well as (2) constructing several spatially local ROMs and coupling them 
to each other.

While the Schwarz alternating method enables the coupling of different non-conformal 
meshes with different resolutions and different time-integration schemes with different
time-steps, as shown in \cite{JLB:mota2017schwarz, JLB:mota2022schwarz}, we 
chose, again for the sake of brevity, to focus the present study on varying the number of POD modes in each ROM subdomain and restricted attention to a non-overlapping DD of $\Omega$ 
into two subdomains, $\Omega_1$ and $\Omega_2$.  As mentioned in the Introduction section,
the present work is the first to formulate and evaluate the non-overlapping variant of our 
Schwarz alternating method for coupling, as all of our past work focused on overlapping coupling \cite{JLB:mota2017schwarz, JLB:mota2022schwarz}.  
For results demonstrating ROM-ROM and FOM-ROM
coupling by means of the overlapping Schwarz alternating method, the interested reader is referred to \cite{JLB:TezaurWCCM}.

Let $\Delta x_i$, $\Delta t_i$, $M_i$ and
$N_{e,i}$
denote the mesh resolution, time-step, POD basis size and number of sample mesh points in subdomain $\Omega_i$ for $i=1,2$, and let $\Delta T$ denote 
the controller time-step.  Recall that the Schwarz tolerance is denoted by $\delta$ (see Section \ref{sec:overlap}).
Our study employed the following DD, mesh-resolutions, time-steps and Schwarz tolerance:
\begin{equation} \label{eq:params}
\begin{array}{c}
    \Omega_1 = [0, 0.6], \hspace{0.2cm} \Omega_2 = [0.6, 1],\\
    \Delta x_1 = \Delta x_2 = 1.0 \times 10^{-3},\\
     \Delta t_1 = \Delta t_2 = \Delta T = 1.0 \times 10^{-7}, \\
     \delta = 1.0\times 10^{-11}.
    \end{array}
\end{equation}
Our choice of DD is motivated by the observation that a much steeper gradient forms in the acceleration
field in the left part of the domain than in the right during the time-interval
considered (see Figure \ref{fig:repro_solns}).  This suggests that it may be possible
to get away with using a relatively small ROM in $\Omega_2$ while still maintaining accuracy,
whereas likely a larger ROM or a FOM will be required in $\Omega_1$. 

In both subdomains, an implicit Newmark-$\beta$ time-integrator with parameters $\beta = 0.49$
and $\gamma = 0.9$ \cite{JLB:Newmark} was used to advance the solution forward
in time\footnote{We chose $\beta = 0.49$
and $\gamma = 0.9$ within the Newmark-$\beta$ scheme, as these parameter values introduce
some numerical dissipation into the discretization, which our FOM-based convergence study revealed
is necessary to stabilize the rapidly-varying acceleration solution.}.
It can be shown \cite{JLB:mota2022schwarz} that the time-step required to resolve the wave and satisfy 
the Courant-Friedrichs-Lewy (CFL) condition is $\Delta t_{\text{CFL}} = \frac{\Delta x}{c}$ where $c = \sqrt{\frac{E}{\rho}}$
is the speed of the wave.  Numerical experiments suggest that the time-step 
required to resolve the wave is in general one order of magnitude smaller than $\Delta t_{\text{CFL}}$.
For the mesh 
resolutions and material parameters being considered, 
$\Delta t_{\text{CFL}} = 1.0\times 10^{-6}$, which justifies our choice of time-step in \eqref{eq:params}.


To assess performance of our coupling method, we have written a MATLAB code that 
implements the discretization methods, model reduction schemes and coupling formulations described above.
In this implementation, the non-negative least squares problem defining the hyper-reduction weights within the ECSW 
algorithm (see Section \ref{sec:ecsw}) is solved using 
MATLAB's {\tt lsqnonneg}
function with an early termination condition, as described in Section \ref{sec:ecsw}.  We present results for both a reproductive (Section \ref{sec:repro}) and a predictive
(Section \ref{sec:pred}) variant of our model problem.  These 
problem variants are described in more detail below.  

To assess
our coupled models' accuracy, we report the mean-square error (MSE) for each of the solution fields: 
displacement, velocity and acceleration.  For the displacement field, 
we calculate this quantity using the following formula: 
\begin{equation} \label{eq:mse_disp}
	\mathcal{E}_{\text{MSE}}(\tilde{~u}_i) := \frac{\sqrt{\sum_{n=1}^S || \tilde{~u}_i^n - ~u_i^n||_2^2}}{\sqrt{\sum_{n=1}^S ||~u_i^n||_2^2}},
\end{equation}
where $~u_i^n$ is the FOM displacement solution (the reference solution) at time $t_n$ in subdomain $\Omega_i$,  and $\tilde{~u}_i^n$ is the ROM 
displacement solution at time $t_n$ in subdomain $\Omega_i$.  
As implied by \eqref{eq:mse_disp},
MSEs for coupled ROM-ROM solutions are calculated by performing 
an analogous coupled FOM-FOM simulation, and computing errors 
in each subdomain ROM with respect to its corresponding
subdomain FOM in the FOM-FOM coupling.  The formula \eqref{eq:mse_disp} 
can also be used to calculate errors in a FOM-ROM coupling with 
respect to a FOM-FOM simulation.  Assuming without loss of generality
that the FOM is prescribed in $\Omega_1$, \eqref{eq:mse_disp} 
is modified by replacing $\tilde{~u}_1$ with the FOM solution in $\Omega_1$
in the FOM-ROM coupling.
The MSE for the velocity
and acceleration fields in subdomain $\Omega_i$ are denoted by $\mathcal{E}_{\text{MSE}}(\tilde{~v}_i)$
and $\mathcal{E}_{\text{MSE}}(\tilde{~a}_i)$, respectively, and defined analogously to \eqref{eq:mse_disp}.

Also reported in the following subsections 
is the total number of Schwarz iterations required for convergence, denoted by $N_{S}$, 
and the CPU time for each run.  All test cases were run in MATLAB in serial on one of two Linux RHEL workstations located at Sandia National Laboratories.  The reproductive cases described in Section \ref{sec:repro} were run on a RHEL8  Intel(R) Xeon(R) CPU E5-2650 v3 \@ 2.30GHz machine, whereas the predictive cases described in Section \ref{sec:pred} were run on a RHEL7 Intel(R) Xeon(R) CPU E7-4880 v2 \@ 2.50GHz machine.
Since
a constant time-step of $1.0 \times 10^{-7}$ was used for 
all of our runs, the average number of Schwarz iterations per time-step can easily be calculated as $N_S/10,000$.

\subsection{Reproductive ROM results}  \label{sec:repro}


We first study the performance of the non-overlapping variant of the proposed 
alternating Schwarz coupling formulation (see Section \ref{sec:nonoverlap}) 
in the reproductive regime.  We prescribe as the initial condition for our problem the Symmetric Gaussian \eqref{eq:ic1}, 
with $a = 1.0\times 10^{-3}$, $b = 0.5$ and $s = 0.02$, and employ the DD, mesh resolutions, time-steps and Schwarz tolerance given in \eqref{eq:params}.
To generate snapshots for building our ROMs, we use the non-overlapping
Schwarz alternating method to perform a FOM-FOM coupled simulation in $\Omega:=\Omega_1 \cup \Omega_2$.  
A total of $S=$10,001 snapshots of the displacement are collected between time 0 and the final time $T = 1.0\times 10^{-3}$
within each subdomain at time-increments of $1.0 \times 10^{-7}$.
These snapshots are used to build POD bases 
for our subdomain ROMs and HROMs.  For comparison and to provide a baseline, we also perform a simulation involving a single-domain
FOM, from which we build a single-domain ROM and HROM. Note that, unless specified otherwise, the snapshots used for POD refer to displacement snapshots.

\begin{figure}[htbp!]
        \begin{center}
      \includegraphics[width=0.7\textwidth]
                      {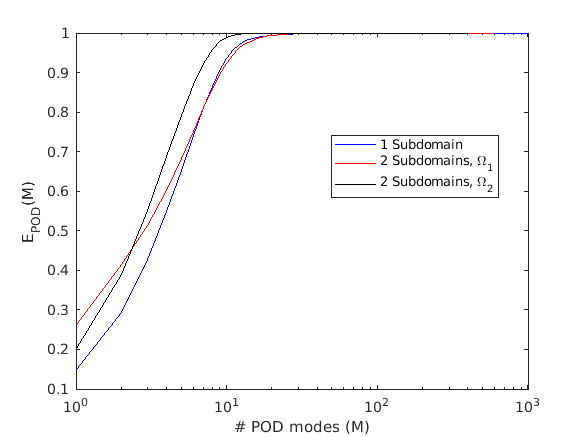}
        \end{center}
        \caption{Energy decay for POD modes constructed from displacement snapshots for the reproductive Symmetric Gaussian variant of the nonlinear wave propagation problem.}
        \label{fig:energy_repro}
\end{figure}

Figure \ref{fig:energy_repro}(a) shows the snapshot energy $\mathcal{E}_{\text{POD}}$ \eqref{eq:energy_pod} 
as a function of the basis size for the single-domain and two subdomain ROM cases.  For the single-subdomain
case, 64 POD modes capture 99.99\% of the snapshot energy.  For the two subdomain
case, 54 and 21 POD modes are required to capture 99.99\% of the snapshot 
energy in $\Omega_1$ and $\Omega_2$, respectively. Interestingly, if one wishes to create a basis 
that captures 100\% of the snapshot energy, significantly more modes are required (509 modes
for the single-domain case, and 414/145 modes in $\Omega_1$/$\Omega_2$ for the two subdomain case).  
As expected, fewer modes are needed to reproduce the solution in $\Omega_2$ than in $\Omega_1$.

Table \ref{tab:reproductive_nonoverlap} reports results for several single-domain ROMs and their hyper-reduced 
variants, termed HROMs, as well as results for several FOM-ROM, ROM-ROM, FOM-HROM and HROM-HROM couplings. 
While the snapshot energies plotted in Figure \ref{fig:energy_repro} suggest that a ROM consisting of less than 70 POD modes should be adequate in reproducing 
the snapshot set, the reader 
can observe by examining Table \ref{tab:reproductive_nonoverlap} that a single-domain ROM comprised of 60 POD modes is fairly inaccurate, achieving an MSE of almost 50\% for the acceleration field.  Since such results are unacceptable in most engineering applications, we decided to consider larger ROMs having $\mathcal{O}(100)$ modes in our coupling studies.

Prior to generating the results summarized in Table \ref{tab:reproductive_nonoverlap}, we performed a parameter sweep study to determine the ``optimal" 
value of the relaxation parameter $\theta$ in \eqref{mechanics_rr2}, that is the value of $\theta$ that yielded the smallest
number of Schwarz iterations, $N_S$, for non-overlapping FOM-FOM and FOM-ROM couplings.  This study 
revealed that the value of $\theta$ that minimizes $N_S$ is 1, which corresponds to \eqref{mechanics_rr2} with no relaxation.
All results reported herein hence used a value of $\theta = 1$. 
The HROMs evaluated were generated using the ECSW hyper-reduction method described in Section \ref{sec:ecsw}.  
A total of $N_h = 20$ snapshots were used to generate the state basis used to train the reduced mesh.  These snapshots
were created by selecting every 500$^{th}$ displacement snapshot used to generate the POD bases.

A number of noteworthy observations can be made by examining Table \ref{tab:reproductive_nonoverlap}.  We begin by discussing the relative performance of the various models.
All couplings summarized in Table \ref{tab:reproductive_nonoverlap} were run on the same Linux RHEL8 workstation, to ensure consistency, as discussed earlier.  
The reader can observe 
by examining this table that all coupled models evaluated 
converged on average in less than 3 Schwarz iterations per time-step.  This indicates
that the coupling method has not introduced a significant amount of overhead. Whereas the FOM-ROM couplings converge in approximately the same number of total Schwarz iterations, $N_S$, as the FOM-FOM coupling, the FOM-HROM, ROM-ROM and HROM-HROM models require more Schwarz iterations to reach convergence. 
It is interesting to observe that the larger (300/80 mode) ROM-ROM is actually slightly faster than the smaller (200/80 mode) ROM-ROM.  While, at first glance, this result may seem counter-intuitive, the reason for this behavior is clear when comparing $N_S$ for the two couplings: the larger ROM-ROM requires fewer Schwarz iterations to converge.  Most likely, this is due to  the larger ROM in $\Omega_1$ being more accurate.
This result is consistent with a previously-observed trend, which showed that coupling less accurate models with each other in general requires more Schwarz iterations to reach 
convergence \cite{JLB:mota2017schwarz, JLB:mota2022schwarz}. 
Despite requiring more Schwarz iterations to converge, the FOM-HROM and HROM-HROM couplings outperform the FOM-FOM coupling in terms of CPU time by between 12.5-32.6\%.  This is a win in the case where a FOM-FOM coupling is the gold standard, as would occur if the problems in $\Omega_1$ and $\Omega_2$ were discretized by two disparate codes, making it impossible to perform a monolithic simulation on  $\Omega$ as a single domain.  While all of the couplings except the smaller HROM-HROM coupling are slower than the single-domain FOM, our preliminary results suggest that speedups over a single-domain FOM are possible through an additive Schwarz implementation of our method (see Remark 3).  
This variant of the Schwarz alternating method will be explored in a subsequent publication.   

\begin{figure}[htbp!]
        \begin{center}
                \subfigure[$t=0$]{
      \includegraphics[width=0.48\textwidth]
                      {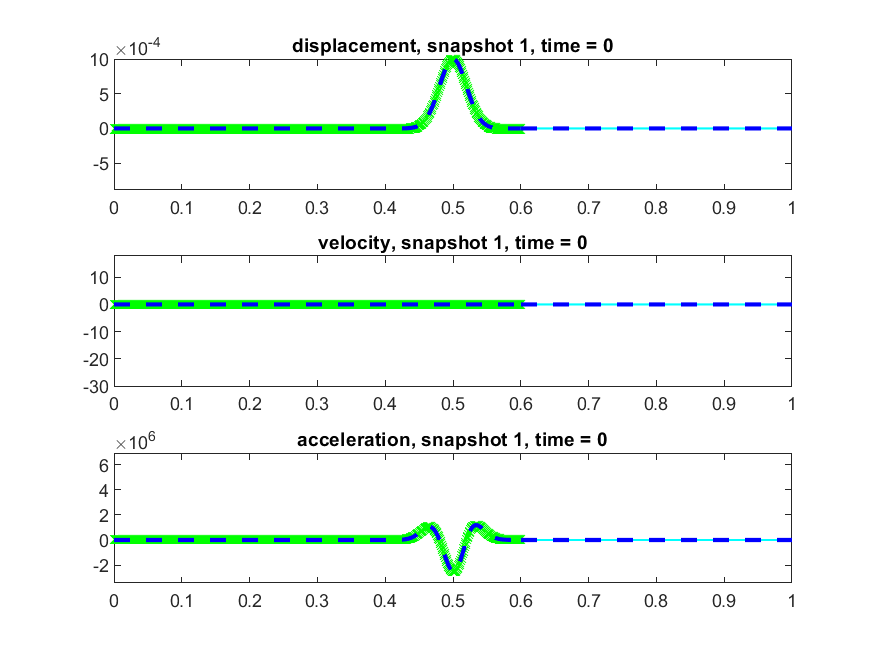}}
	      \subfigure[$t=2.5\times 10^{-4}$]{
      \includegraphics[width=0.48\textwidth]
                      {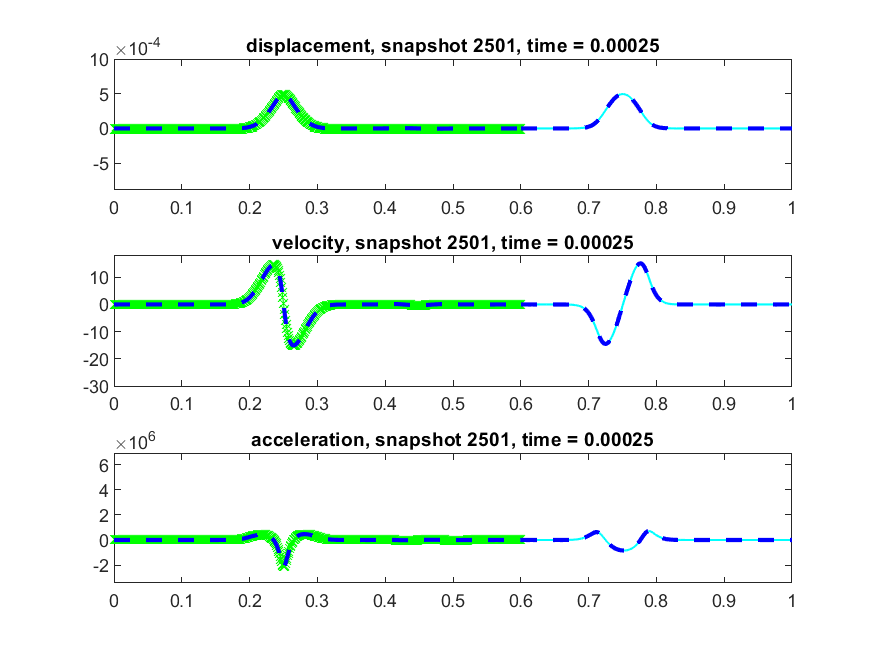}}
	      \subfigure[$t=7.5\times 10^{-4}$]{
      \includegraphics[width=0.48\textwidth]
                      {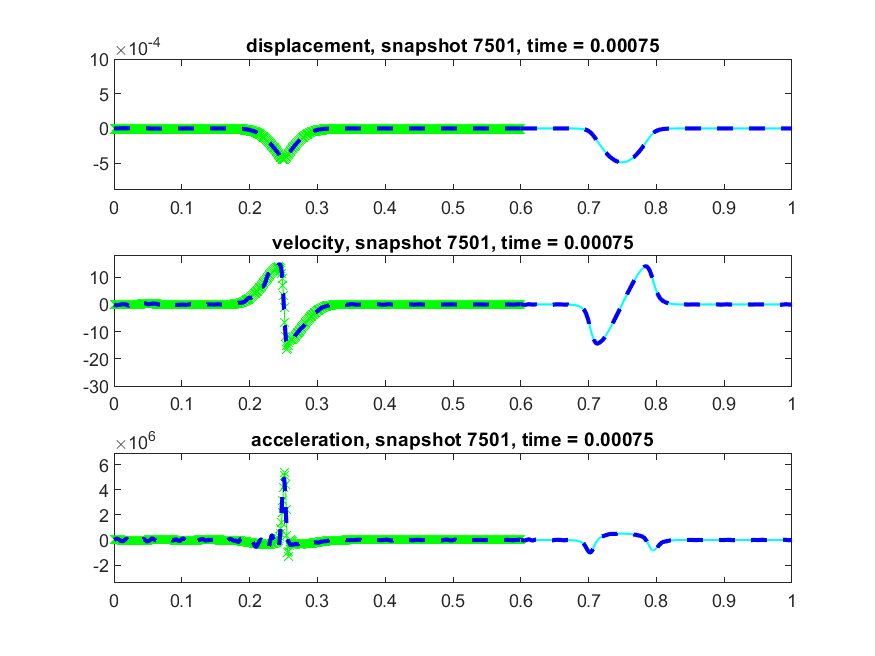}}
	      \subfigure[$t=1.0 \times 10^{-3}$]{
      \includegraphics[width=0.48\textwidth]
                      {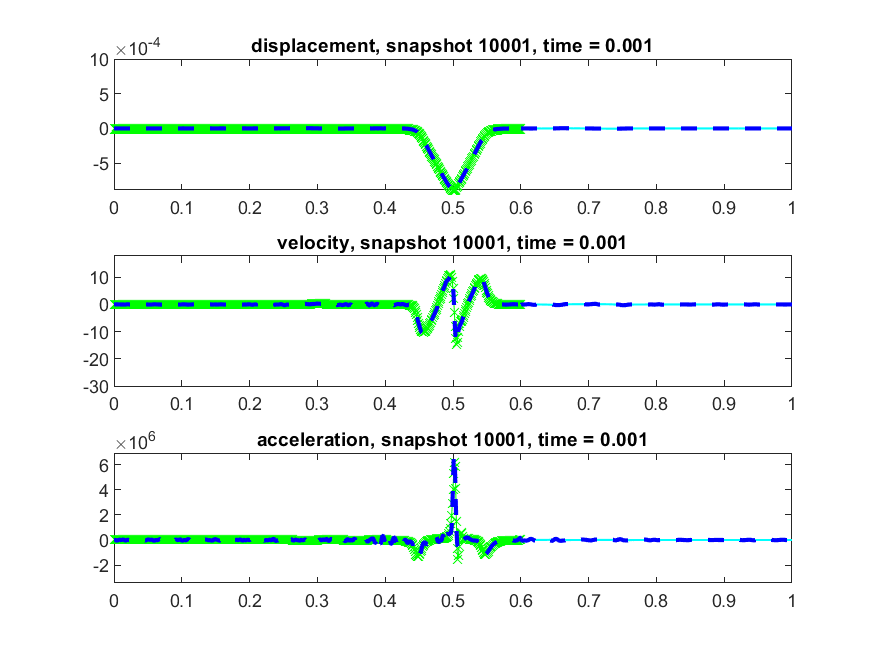}}
        \end{center}
	\caption{Plots of the HROM-HROM solutions with $M_1 = 200$, $M_2 = 80$, $N_{e,1}=315$ and $N_{e,2} = 130$ (see Table \ref{tab:reproductive_nonoverlap}) compared to a single-domain
	FOM solution.  The coupled solutions in $\Omega_1$ and $\Omega_2$ are shown in green and cyan, respectively,
	whereas the single-domain FOM solution is shown in blue.  No coupling artifacts are observed in the displacement, velocity and acceleration solutions.}
        \label{fig:repro_solns}
\end{figure}

Having discussed performance, we now turn our attention to accuracy.  First, convergence is observed with basis refinement for all the models except the larger HROM-HROM.  The reader can observe that the MSEs in $\Omega_2$ are lower for all ROM-ROM and HROM-HROM couplings despite the fact that fewer modes are employed within this subdomain.  This is due to the solution in $\Omega_2$ being far smoother than the solution in $\Omega_1$ for the time interval considered (see Figure \ref{fig:repro_solns}).  Among the most accurate models involving HROMs are the FOM-HROMs considered, which achieve errors of $\mathcal{O}(0.01\%)$ or less in the displacement solution and of  $\mathcal{O}(1\%)$ or less in the acceleration solution.   It is interesting and encouraging to remark that the larger FOM-ROM model is incredibly accurate, achieving an MSE of $\mathcal{O}(10^{-11})$ in the displacement solution.   While our HROM-HROM couplings achieve reasonable errors for the displacement and velocity fields, non-trivial errors of $\mathcal{O}(10\%)$ appear to be unavoidable when using these models given the fixed early termination criterion used in our implementation of ECSW.  In interpreting 
this result, it is important to recognize that all couplings involving ROMs and HROMs are at least as accurate as the single-domain ROMs and HROMs evaluated in Table \ref{tab:reproductive_nonoverlap}.  This result suggests that the errors we are seeing in the coupled models are due to the inaccuracy of the individual subdomain models, rather than errors introduced by our coupling framework. 
It may be possible to improve the accuracy of coupled ROMs/HROMs by decomposing the spatial domain into more subdomains and creating/coupling a larger number of spatially-localized ROMs/HROMs.  
It may also be possible to improve the HROM-HROM results in Table \ref{tab:reproductive_nonoverlap} by simply changing the termination tolerance  within MATLAB's {\tt lsqnonneg} algorithm used to solve \eqref{eq:ecsw_solve}.  As discussed in Section \ref{sec:ecsw}, we did not do this here, as our goal was to have a consistent termination criterion for all HROMs being evaluated.

\begin{table}[h!]
  \begin{center}
    \caption{Reproductive non-overlapping Schwarz coupling results, with $\theta = 1$.  
	  The POD modes
	  for the couplings involving ROMs and HROMs were constructed from snapshots of only the displacement field.  Couplings outperforming the FOM-FOM model in terms of CPU-time 
	  and with reasonable errors are high-lighted in green.
	  }
    \label{tab:reproductive_nonoverlap}
    \begin{scriptsize}
    \begin{tabular}{c|c|c|c|c|c|c|c}
	    \multirow{2}{*}{Model} & \multirow{2}{*}{$M_1$/$M_2$} & \multirow{2}{*}{$N_{e,1}$/$N_{e,2}$} & CPU  & $\mathcal{E}_{\text{MSE}}(\tilde{~u}_1)$/ & $\mathcal{E}_{\text{MSE}}(\tilde{~v}_1)$/ & $\mathcal{E}_{\text{MSE}}(\tilde{~a}_1)$/ & \multirow{2}{*}{$N_S$} \\
	    & & & time (s) & $\mathcal{E}_{\text{MSE}}(\tilde{~u}_2)$ & $\mathcal{E}_{\text{MSE}}(\tilde{~v}_2)$ & $\mathcal{E}_{\text{MSE}}(\tilde{~a}_2)$\\
      \hline
	    \hline
	    FOM & $-$/$-$ & $-$/$-$ & $1.871\times 10^3$& $-$/$-$ & $-$/$-$ & $-$/$-$ & $-$ \\
	    \hline 
	    ROM & 60/$-$ & $-$/$-$ & $1.398\times 10^3$&$1.659\times 10^{-2}$/$-$ &$1.037\times 10^{-1}$ & $4.681\times 10^{-1}$/$-$ & $-$\\
	    \hline
	      HROM & 60/$-$ & 155/$-$ & $5.878\times 10^2$ & $1.730\times 10^{-2}$/$-$ &$1.063\times 10^{-1}$/$-$ & $4.741\times 10^{-1}$/$-$ & $-$\\
	      \hline
	    ROM & 200/$-$ & $-$/$-$ &$1.448\times 10^3$ & $2.287\times 10^{-4}/-$&$4.038\times 10^{-3}$/$-$ &$4.542\times 10^{-2}$/$-$  & $-$\\
	    \hline
	    HROM & 200/$-$ & 428/$-$ &$9.229\times 10^2$ &$8.396\times 10^{-4}$/$-$ &$8.947\times 10^{-3}$/$-$ &$7.462\times 10^{-2}$/$-$  & $-$\\
	    \hline \hline

  FOM-FOM & $-$/$-$ &$-$/$-$ & $2.345\times 10^{3}$ & $-$ & $-$ & $-$ & 24,630 \\
	    \hline
	    \cellcolor{green!25}FOM-ROM & \cellcolor{green!25}$-$/80 & 
	    \cellcolor{green!25}$-$/$-$ & \cellcolor{green!25} $2.341\times 10^3$ & \cellcolor{green!25}$2.171\times 10^{-6}$/ & \cellcolor{green!25}$3.884\times 10^{-5}$/ & \cellcolor{green!25}$2.982\times 10^{-4}$/ & \cellcolor{green!25} 25,227 \\
	    \cellcolor{green!25}& \cellcolor{green!25}& \cellcolor{green!25}& \cellcolor{green!25}& \cellcolor{green!25}$1.253\times 10^{-5}$& \cellcolor{green!25}$2.401\times 10^{-4}$ & \cellcolor{green!25}$2.805\times 10^{-3}$ & \cellcolor{green!25}\\
	    \hline
	    \cellcolor{green!25}FOM-HROM & \cellcolor{green!25}$-$/80 & \cellcolor{green!25}$-$/130 & \cellcolor{green!25}$2.085\times 10^3$ & \cellcolor{green!25}$2.022\times 10^{-4}$/ &\cellcolor{green!25} $1.723e\times 10^{-3}$/ & \cellcolor{green!25}$7.421\times 10^{-3}$/ & \cellcolor{green!25}29,678 \\
	    \cellcolor{green!25}& \cellcolor{green!25}& \cellcolor{green!25}& \cellcolor{green!25}&\cellcolor{green!25} $5.734\times 10^{-4}$& \cellcolor{green!25}$5.776\times 10^{-3}$ &\cellcolor{green!25} $3.791\times 10^{-2}$ & \cellcolor{green!25}\\
	    \hline
	    \multirow{2}{*}{FOM-ROM} & \multirow{2}{*}{$-$/200} & \multirow{2}{*}{$-$/$-$} & \multirow{2}{*}{$2.449\times 10^3$} & $4.754\times 10^{-12}$/ & $1.835\times 10^{-10}$/ & $5.550\times 10^{-9}$/ & \multirow{2}{*}{24,630} \\
	    & & & & $7.357\times 10^{-11}$ & $4.027\times 10^{-9}$ & $1.401\times 10^{-7}$\\
	    \hline
	    \multirow{2}{*}{FOM-HROM} & \multirow{2}{*}{$-$/200} & \multirow{2}{*}{$-$/252} & \multirow{2}{*}{$2.352\times 10^3$} & $1.421\times 10^{-5}$/ & $1.724\times 10^{-4}$/ & $9.567\times 10^{-4}$/ & \multirow{2}{*}{27,156}  \\
	    & & & & $4.563\times 10^{-4}$& $2.243\times 10^{-3}$ & $1.364\times 10^{-2}$\\
	    \hline
	    \multirow{2}{*}{ROM-ROM} & \multirow{2}{*}{200/80} & \multirow{2}{*}{$-$/$-$} & \multirow{2}{*}{$2.778\times 10^3$} & $4.861\times 10^{-5}$/ & $1.219\times 10^{-3}$/ & $1.586\times 10^{-2}$/ & \multirow{2}{*}{27,810}  \\
	    & & & & $3.093\times 10^{-5}$& $4.177\times 10^{-4}$ & $3.936\times 10^{-3}$\\
	    \hline
	    \cellcolor{green!25}HROM-HROM & \cellcolor{green!25}200/80 & \cellcolor{green!25}315/130 & \cellcolor{green!25}$1.769\times 10^3$ & \cellcolor{green!25}$3.410\times 10^{-3}$/ &\cellcolor{green!25}$4.110\times 10^{-2}$/ & \cellcolor{green!25}$2.485\times 10^{-1}$/ & \cellcolor{green!25}29,860  \\
	    \cellcolor{green!25}& \cellcolor{green!25}& \cellcolor{green!25}& \cellcolor{green!25}& \cellcolor{green!25}$6.662\times 10^{-4}$ & \cellcolor{green!25}$6.432\times 10^{-3}$ & \cellcolor{green!25}$4.307\times 10^{-2}$ & \cellcolor{green!25}\\
	    \hline
	    \multirow{2}{*}{ROM-ROM} & \multirow{2}{*}{300/80} & \multirow{2}{*}{$-$/$-$} & \multirow{2}{*}{$2.646\times 10^3$} & $2.580\times 10^{-6}$/ & $6.226\times 10^{-5}$/ & $9.470\times 10^{-4}$ & \multirow{2}{*}{25,059} \\
	    & & & & $1.292\times 10^{-5}$ & $2.483\times 10^{-4}$ & $2.906\times 10^{-3}$\\
	    \hline

	    \cellcolor{green!25} HROM-HROM & \cellcolor{green!25}300/80 & \cellcolor{green!25}405/130 & \cellcolor{green!25}$1.938\times 10^3$ & \cellcolor{green!25}$6.960\times 10^{-3}$ & \cellcolor{green!25}$6.328\times 10^{-2}$ & \cellcolor{green!25}$3.137\times 10^{-1}$ &\cellcolor{green!25}29,896  \\
	   \cellcolor{green!25} &\cellcolor{green!25} & \cellcolor{green!25}& \cellcolor{green!25}& \cellcolor{green!25} $7.230\times 10^{-4}$ & \cellcolor{green!25}$7.403\times 10^{-3}$ & \cellcolor{green!25}$4.960\times 10^{-2}$ & \cellcolor{green!25}

    \end{tabular}
    \end{scriptsize}
  \end{center}
\end{table}

\begin{figure}[htbp!]
        \begin{center}
      \includegraphics[width=1\textwidth]
             {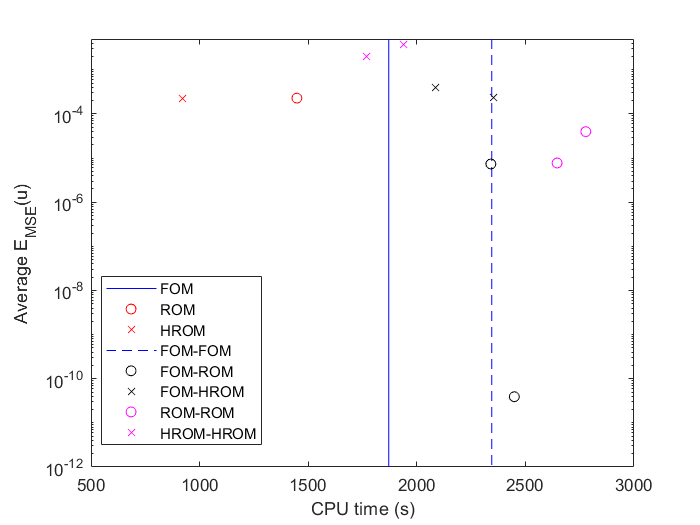}
        \end{center}
	\caption{Pareto plot showing CPU time in seconds versus the average displacement MSE over all subdomains being coupled for the reproductive couplings evaluated in Table \ref{tab:reproductive_nonoverlap}.}
        \label{fig:repro_pareto}
\end{figure}

In an effort to provide a complete picture of the relative computational cost and accuracy of the models evaluated, we show in Figure \ref{fig:repro_pareto} a Pareto plot, which plots the CPU time in seconds versus the average displacement MSE over all subdomains being coupled.  The reader can observe that, while the single-domain ROMs and HROMs are the fastest, our coupling methodology enables us to achieve lower errors by performing ROM-ROM and FOM-ROM couplings. Future work will examine ways to improve the accuracy and efficiency of FOM-HROM and HROM-HROM couplings, which are not optimal according 
to Figure \ref{fig:repro_pareto}.  It is not possible to make a definitive conclusion about the general utility of FOM-HROM and HROM-HROM couplings  until we have evaluated our coupling methodology in two or three spatial dimensions.  In multiple spatial dimensions, we anticipate seeing a greater benefit from hyper-reduction, as well as a greater potential for generating  ``optimal" (MSE and CPU-time minimizing) DDs and ROM/FOM assignments, in the spirit of \cite{JLB:Bergmann2018}.  

In Figure \ref{fig:repro_solns}, we plot the solution 
computed using our non-overlapping Schwarz coupling method for the smaller HROM-HROM coupling summarized
in Table \ref{tab:reproductive_nonoverlap}.  In this figure, the solution in $\Omega_1$ is shown in green, 
and the solution in $\Omega_2$ is shown in cyan.  It is evident that the coupled 
solutions are indistinguishable from the single-domain FOM solution in the full domain $\Omega$, 
shown in blue.  This indicates  that the coupling method has not introduced any spurious artifacts
into the discretization.

\subsubsection{Alternate POD bases} \label{sec:disp-acce-POD}

The reader can observe by inspecting Figure \ref{fig:repro_solns} that, whereas the displacement and velocity solutions are relatively smooth, 
sharp gradients form and propagate in the acceleration field over the course of the simulation.  Since these gradients are difficult
to capture using POD modes, significantly larger MSEs are observed in the acceleration field than in the displacement 
and velocity fields (see Table \ref{tab:reproductive_nonoverlap}).  This observation motivated us to explore a simple idea for improving 
on the results in Table \ref{tab:reproductive_nonoverlap} by augmenting the POD basis used in 
constructing each of our ROMs with modes calculated using snapshots of the acceleration field.  Let $~\Phi_M^{~u}$ and $~\Phi_M^{~a}$
denote POD bases of size $M$ computed from displacement and acceleration modes, respectively.  Our procedure 
for generating a combined displacement-acceleration POD bases involved first computing  $~\Phi_M^{~u}$ and $~\Phi_M^{~a}$ independently.
Once this was done, an SVD was performed of the augmented matrix $\left[~\Phi_M^{~u}, ~\Phi_M^{~a}\right]$ to remove potential linear 
dependencies between the two bases, and the resulting basis was truncated as desired to yield a displacement-acceleration POD 
basis.  \\

\noindent {\bf \textit{Remark 4.}} The idea to include snapshots of the acceleration field
within the POD basis of a ROM is related to the idea of including so-called ``time difference quotients" (numerical estimates 
of the time-derivative of the primary solution field) in the generation of POD modes 
\cite{JLB:Iliescu2014}. These difference quotients can be thought of as representing the right-hand side of a nonlinear PDE being solved, and 
were demonstrated in \cite{JLB:Iliescu2014} to lead to ROMs with improved convergence rates.\\

Unfortunately, while the ``displacement-acceleration" basis generation approach discussed above improved the errors by approximately four orders of magnitude for the FOM-ROM couplings we tried, the approach was actually deleterious when it came to FOM-HROM, ROM-ROM and HROM-HROM couplings.  The reason for this is unknown at the present time, and is something we plan to look into in future work.

\subsection{Predictive ROM results} \label{sec:pred} 

Having evaluated the proposed coupling methodology on a reproductive test case, we now turn our attention to a 
predictive variant of the nonlinear wave propagation problem.
First, snapshots are generated by simulating the problem with the Symmetric Gaussian \eqref{eq:ic1} initial
condition with parameters $a = 1.0\times 10^{-3}$, $b = 0.5$ and $s = 0.02$ up to time $T = 1.0\times 10^{-3}$, as before.  A total of 10,001 snapshots of the displacement field are generated at increments of $1.0\times 10^{-7}$.
POD modes are constructed from these snapshots, and used to predict 
the solution of our model problem with a different initial condition, namely the Rounded Square \eqref{eq:ic2}, with parameters $a = 5.0\times 10^{-4}$, $b = 100$ and $s = 0.6$.  As with the Symmetric Gaussian 
variant of this problem, the predictive simulation is 
run until time $T = 1.0\times 10^{-3}$.

\begin{figure}[htbp!]
        \begin{center}
      \includegraphics[width=0.7\textwidth]
                      {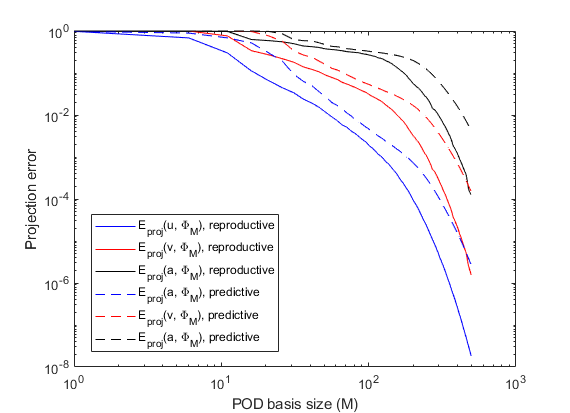}
        \end{center}
        \caption{Projection errors for the nonlinear wave propagation problem with the Rounded Square initial condition.}
        \label{fig:proj_errs}
\end{figure}

Before building any ROMs and performing any couplings, we assess the projection error for the Rounded Square 
nonlinear wave propagation problem, defined as
\begin{equation} \label{eq:proj_err_disp}
 \mathcal{E}_{\text{proj}}(~u, ~\Phi_M) := \frac{||~u - ~\Phi_M(~\Phi_M^T ~\Phi_M)^{-1} ~\Phi_M^T ~u||_2}{||~u||_2},   
\end{equation}
where $~u$ denotes the snapshots of the displacement field.  The projection error can be 
computed in a similar way for the velocity and acceleration snapshots, and is denoted by $\mathcal{E}_{\text{proj}}(~v, ~\Phi_M)$
and $\mathcal{E}_{\text{proj}}(~a, ~\Phi_M)$, respectively.  
Equation \eqref{eq:proj_err_disp} provides a straightforward and inexpensive way to evaluate a given basis's ability to represent a solution without having to construct and run an entire projection-based ROM.   
Figure \ref{fig:proj_errs} shows a plot of the projection error 
for the Rounded Square nonlinear wave propagation problem simulated on one domain $\Omega$ as a function 
of the basis size $M$.  We calculate and report $\mathcal{E}_{\text{proj}}(~u, ~\Phi_M)$, $\mathcal{E}_{\text{proj}}(~v, ~\Phi_M)$ and $\mathcal{E}_{\text{proj}}(~a, ~\Phi_M)$
for two kinds of POD bases: reproductive and predictive.  The reader can observe that the projection error decreases with
basis refinement for both basis types, indicating that obtaining a reasonably accurate solution should be possible even
with the predictive POD basis, provided enough modes are retained.  It is noted that the acceleration projection error
is roughly two order of magnitude higher than the displacement projection error, and exhibits a much slower decay 
than the displacement and velocity projection errors.

As for the reproductive problem considered in Section \ref{sec:repro}, we evaluate the performance 
of the non-overlapping variant of the proposed Schwarz-based coupling method (Section \ref{sec:nonoverlap}) with the DDs, mesh resolutions, time-steps and Schwarz tolerance
given in \eqref{eq:params}.
Again, we set $\theta = 1$ in \eqref{eq:relaxation}, as this value yielded the best performance.
We consider predictive ROMs with a relatively large number of POD modes in each subdomain: 300 modes for 
a single-domain ROM, and 300/200 modes for a two subdomain ROM in $\Omega_1/\Omega_2$.  This 
choice of basis size is motivated by the projection error results plotted in Figure \ref{fig:proj_errs}, 
which demonstrate that a minimum of $O(100)$ modes are needed to reproduce our three solution fields to a sufficient
accuracy in the predictive regime.    

The main results are summarized in Table \ref{tab:predictive}.  All runs summarized in Table \ref{tab:predictive} were performed on the same Linux RHEL7 workstation, to ensure consistency of CPU times for an apples-to-apples comparison, as described earlier.  First, we assess the single-domain ROM and HROM results.  The reader can observe that the MSEs for the ROM and HROM are almost identical.  This is likely due to the fact that a very large number of mesh points (more than 60\%) is being selected by the ECSW algorithm when constructing the HROM.  It is interesting to remark that the ROM is actually slower than the FOM.  This happens frequently when constructing projection-based ROMs for nonlinear problems without using hyper-reduction, and highlights the importance of hyper-reduction.

Turning our attention to the coupled results, we begin by discussing accuracy, as done in Section \ref{sec:repro}.   It can be 
seen that the FOM-ROM coupling is remarkably accurate, achieving MSEs as low as $\mathcal{O}(10^{-8})$ for this \textit{predictive} problem.  This result supports our claim that a more accurate and robust model can be obtained via our Schwarz-based coupling, and suggests that no coupling errors have been introduced into the discretization.  Although the FOM-HROM and ROM-ROM 
couplings are less accurate than the FOM-ROM coupling, one can 
see that they achieve MSEs lower than the single-domain ROM and single-domain HROM summarized in Table \ref{tab:predictive}\footnote{There is, of course, the caveat 
that the single-domain ROM has less total modes than our ROM-ROM coupling being evaluated.}.  The MSEs for the HROM-HROM model in Table \ref{tab:predictive} are about one order of magnitude higher than its corresponding ROM-ROM, but nonetheless on par with the single-domain HROM we evaluated when it comes to MSE.

Next, we discuss the efficiency of the various models in Table \ref{tab:predictive}.
Interestingly, the FOM-ROM model takes slightly less CPU time to converge than its corresponding FOM-FOM model despite not having hyper-reduction, and requires approximately the same number of Schwarz iterations to converge.  The same cannot be said for the FOM-HROM, ROM-ROM and HROM-HROM summarized in Table \ref{tab:predictive}, which are slower than their analogous FOM-FOM coupled model.  This is due to the number of Schwarz iterations, which is up to 30\% higher than for the FOM-FOM and FOM-ROM couplings (with an average of as many as three Schwarz iterations per time-step needed to converge the HROM-HROM solution).  As discussed in Section \ref{sec:repro}, the larger number of Schwarz iterations needed to reach convergence is most likely due to the lower accuracy of the individual models being coupled.  In the case of the HROM-HROM case, the reader can observe that, like for the single-domain HROM, more than 60\% of the mesh points are being sampled by the ECSW algorithm for the HROM-HROM coupling, which contributes to the higher CPU time attained by this model.  
The reader is reminded that, although the non-negative least squares solver \texttt{lsqnonneg} in MATLAB used in the ECSW procedure has a consistent early termination condition given by a fixed solution step size tolerance of $10^{-4}$, the relative error tolerance $||~C ~\xi - ~d ||_2 / ||~d||_2$ from Section \ref{sec:ecsw} will differ between reduced meshes depending on the choice of number of modes and number of degrees of freedom in a given subdomain.  It may be possible to obtain a more accurate HROM solution by simply tweaking the {\tt lsqnonneg} termination tolerance, but we chose not to do that here to ensure consistency between all methods being evaluated, as discussed in Section \ref{sec:ecsw}. 
Like in Section \ref{sec:repro}, we remark that all coupled models in Table \ref{tab:predictive} can be improved further
through the use of the additive Schwarz variant, as discussed in Remark 3.

\begin{table}[h!]
  \begin{center}
    \caption{Predictive non-overlapping Schwarz coupling results, with $\theta = 1$.  
	  The POD modes
	  for the couplings involving ROMs and HROMs were constructed from snapshots of only the displacement field.  
	  The ROMs and HROMs summarized in this table had the following numbers of POD modes: $M_1 = 300$, $M_2 = 200$.  Couplings outperforming the FOM-FOM model in terms of CPU-time 
	  and with reasonable errors are high-lighted in green.}
    \label{tab:predictive}
     \begin{scriptsize}
    \begin{tabular}{c|c|c|c|c|c|c}
	    \multirow{2}{*}{Model} & CPU  & \multirow{2}{*}{$N_{e, 1} / N_{e, 2}$} & $\mathcal{E}_{\text{MSE}}(\tilde{~u}_1)$/ & $\mathcal{E}_{\text{MSE}}(\tilde{~v}_1)$/ & $\mathcal{E}_{\text{MSE}}(\tilde{~a}_1)$/ & \multirow{2}{*}{$N_S$} \\
	    & time (s) & & $\mathcal{E}_{\text{MSE}}(\tilde{~u}_2)$ & $\mathcal{E}_{\text{MSE}}(\tilde{~v}_2)$ & $\mathcal{E}_{\text{MSE}}(\tilde{~a}_2)$\\
      \hline
	    \hline
	   FOM& 1.288 $\times 10^3$ & $-$/$-$ & $-$/$-$ & $-$/$-$ & $-$/$-$ & $-$  \\
	   \hline
	   ROM & 1.358 $\times 10^3$ & $-$/$-$ & 3.451$\times 10^{-3}$/$-$ & $6.750 \times 10^{-2}$/$-$ & $3.021 \times 10^{-1}$/$-$ & $-$\\
	   \hline
	   HROM & $9.759 \times 10^2$ & $614$/$-$ & $3.463 \times 10^{-3}$/$-$ & $6.750 \times 10^{-2}$/$-$ & $3.021 \times 10^{-1}$/$-$ & $-$\\
	   \hline \hline
	   FOM-FOM  & $2.133 \times 10^3$ & $-$/$-$ & $-$/$-$ & $-$/$-$ & $-$/$-$& 23,280   \\
	   \hline
	   \cellcolor{green!25} FOM-ROM &  \cellcolor{green!25}$2.084 \times 10^3$ &  \cellcolor{green!25}$-/-$ &  \cellcolor{green!25}$1.907\times 10^{-8}$/ &  \cellcolor{green!25}$1.461 \times 10^{-6}$/ &  \cellcolor{green!25}$3.973\times 10^{-5}$/ &  \cellcolor{green!25}23,288 \\
	    \cellcolor{green!25} & \cellcolor{green!25} & \cellcolor{green!25} &  \cellcolor{green!25}$1.170\times 10^{-6}$ &  \cellcolor{green!25}$9.882\times 10^{-5}$ &  \cellcolor{green!25}$1.757\times 10^{-3}$ &  \cellcolor{green!25} \\
	    \hline
	   \multirow{2}{*}{FOM-HROM} & \multirow{2}{*}{$2.219 \times 10^3$} & \multirow{2}{*}{$-/253$} & 1.967 $\times 10^{-4}$ & $4.986 \times 10^{-3}$ & $2.768 \times 10^{-2}$ & \multirow{2}{*}{29,700}\\ 
       & & & $1.720 \times 10^{-3}$ & 4.185 $\times 10^{-2}$ & $2.388 \times 10^{-1}$\\
	   \hline 
	   \multirow{2}{*}{ROM-ROM} & \multirow{2}{*}{$2.502 \times 10^3$} & \multirow{2}{*}{$-/-$} & $5.592 \times 10^{-4}$/ & $1.575\times 10^{-2}$/ & $9.197\times 10^{-2}$/ & \multirow{2}{*}{26,220} \\
	    & & & $4.346\times 10^{-4}$ & $1.001\times 10^{-2}$ & $5.304\times 10^{-2}$\\
	   \hline
	   \multirow{2}{*}{HROM-HROM}  & \multirow{2}{*}{2.200$\times 10^3$} & \multirow{2}{*}{$405/253$} & $4.802\times 10^{-3}$ & $8.500\times 10^{-2}$ & $3.744\times 10^{-1}$ & \multirow{2}{*}{30,067}\\ 
    & & & $1.960 \times 10^{-3}$ & 4.630 $\times 10^{-2}$ & 2.580 $\times 10^{-1}$
	   \\
    \end{tabular}
    \end{scriptsize}
  \end{center}
\end{table}

\begin{figure}[htbp!]
        \begin{center}
      \includegraphics[width=1\textwidth]
             {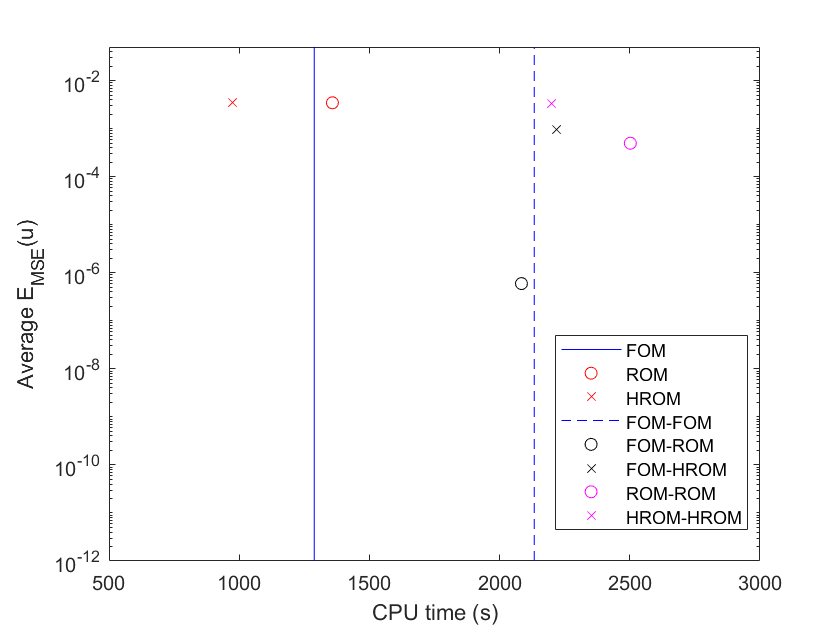}
        \end{center}
	\caption{Pareto plot showing CPU time in seconds versus the average displacement MSE over all subdomains being coupled for the predictive couplings evaluated in Table \ref{tab:predictive}.}
        \label{fig:pred_pareto}
\end{figure}

As for our reproductive test case, we combine the accuracy and efficiency results summarized in Table \ref{tab:predictive} into a Pareto plot, which shows the CPU time in seconds versus the average displacement MSE over all subdomains being coupled (Figure \ref{fig:pred_pareto}).  It is clear from this figure that, by coupling a ROM to a FOM, it is possible to reduce error of a given reduced model by several orders of magnitude.  Figure \ref{fig:pred_pareto} reinforces the need to investigate ways to improve FOM-HROM and HROM-HROM efficiency and accuracy in future work.

\begin{figure}[htbp!]
        \begin{center}
                \subfigure[Predictive single-domain ROM]{
      \includegraphics[width=0.96\textwidth]
                      {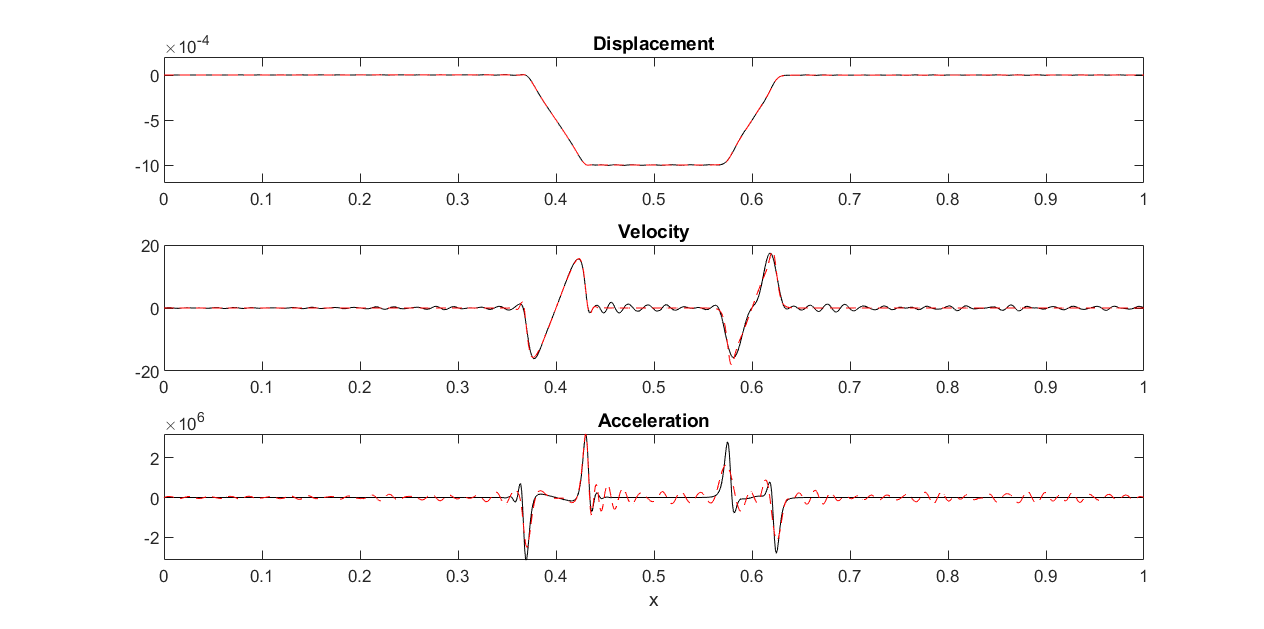}}
                           \subfigure[Predictive FOM-HROM]{
      \includegraphics[width=0.96\textwidth]
                      {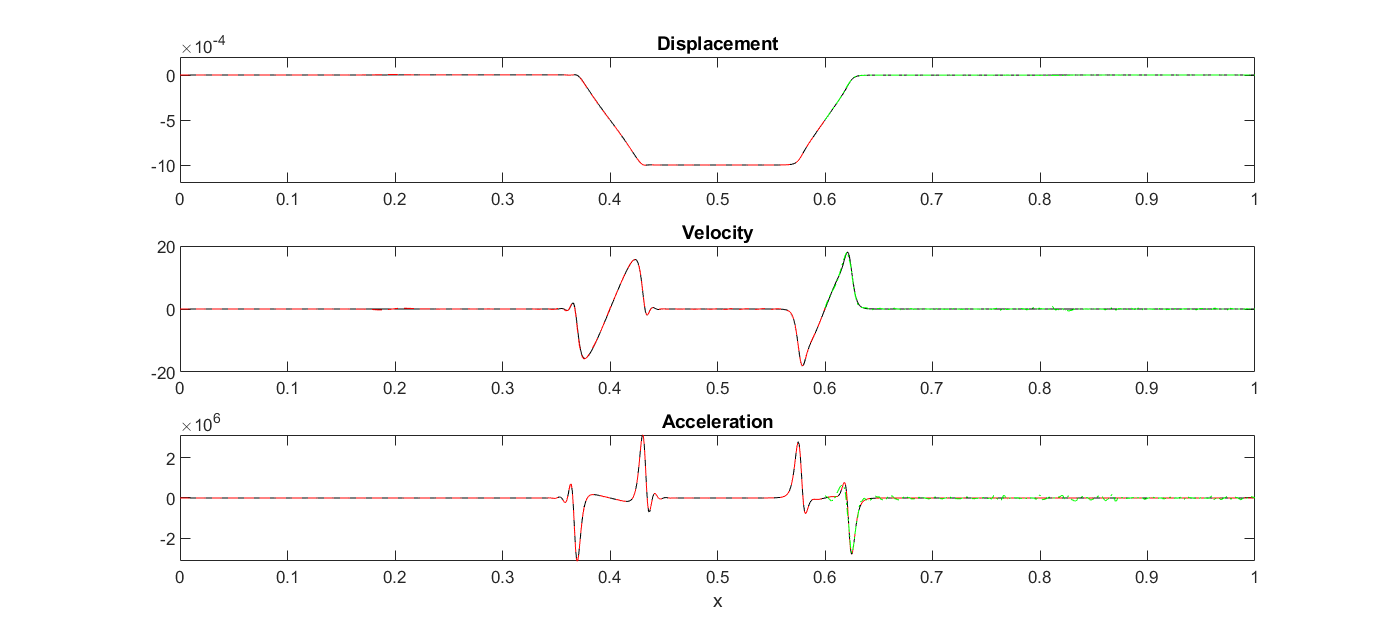}}
        \end{center}
	\caption{Solutions to the predictive Rounded Square nonlinear wave propagation problem at the final time $T = 1.0\times 10^{-3}$.  The coupled solutions in $\Omega_1$ and $\Omega_2$ are shown in red and green, respectively; the single-domain FOM solution is shown in black. Whereas the velocity and acceleration components of the single-domain ROM is fraught with oscillations, the FOM-HROM coupled solutions agrees closely with the single-domain FOM solution.}
        \label{fig:pred_solns}
\end{figure}

Figure \ref{fig:pred_solns} shows the displacement, velocity and acceleration solutions 
at the final time $T = 1.0\times 10^{-3}$ for several of the predictive models being evaluated: the single-domain ROM and the coupled FOM-HROM. The coupled solutions in $\Omega_1$ and $\Omega_2$ are shown in red and green, respectively.  These are plotted on top of a single-domain FOM solution, shown in black.  
The reader can observe that the predictive single-domain ROM solution (Figure \ref{fig:pred_solns}(a)) exhibits some spurious oscillations 
in the velocity and acceleration fields. In contrast, the FOM-HROM solution (Figure \ref{fig:pred_solns}(b)) is smooth and in good agreement with the single-domain FOM.

While the results summarized in Table \ref{tab:predictive} and Figure \ref{fig:pred_solns} are promising, some work still needs to be done in improving the accuracy of the coupled model when the Schwarz alternating method is used to stitch together HROMs.  We will explore strategies to do this in the context of a two-dimensional (2D) problem in a subsequent publication.

\section{Summary and future work} \label{sec:conc}

In this paper, we described a methodology for coupling projection-based ROMs with each other and with conventional high-fidelity finite element models by means of the Schwarz alternating method. In this method, the physical domain is decomposed into two or more subdomains, and a sequence of subdomain-local problems is solved, with information propagating through carefully-constructed transmission boundary conditions posed on the subdomain boundaries.  
The method was formulated for both overlapping and non-overlapping domain decompositions in the context of a generic nonlinear solid dynamics problem.  While our numerical experiments focused on couplings in which the same time-integrator and time-step is used in all subdomains being coupled, we emphasize that the method is capable of coupling not only disparate discretizations but also different time-integrators with disparate time-scales, as demonstrated in \cite{JLB:mota2022schwarz, JLB:TezaurWCCM}.

The utility of the proposed coupling approach is two-fold.  First, it provides a mechanism for enabling the ``plug-and-play" integration of data-driven models  into existing multi-scale and multi-physics coupling frameworks, with minimal intrusion (i.e., through the introduction of a simple outer loop around the two or more disparate models being coupled).  Second,  couplings such as those performed herein have the potential of improving the predictive viability of projection-based ROMs, by enabling the spatial localization of ROMs (via domain decomposition) and the online integration of high-fidelity information into these models (via FOM coupling). 

Our numerical results 
for a 1D nonlinear wave propagation problem (a problem which poses a number of challenges for traditional POD-based model reduction approaches) are promising.  These results demonstrate that the proposed methodology is capable of coupling disparate models without introducing numerical artifacts into the solution for both reproductive and predictive problems.  Additionally, they show that the resulting couplings can be cost-effective when combined with hyper-reduction.   Our conclusion that  FOM-ROM couplings are particularly accurate 
for the problems considered is not surprising, given the sharp gradients present in the acceleration component of the solution.  Features such as these are incredibly 
difficult to capture using POD modes alone, especially when they propagate dynamically in time.  Our results suggest that some more 
work can be done in the future to try to improve the accuracy and efficiency of couplings involving HROMs.  This task will be pursued in the context of 2D and three-dimensional (3D) problems with localized features in their solutions (e.g., shedding vortices behind a cylinder in a fluid simulation, failure and strain localization in a solid simulation), where more benefits are expected from both hyper-reduction as well as domain decomposition.
Although attention herein was restricted to non-overlapping couplings, 
similar results were obtained with the overlapping version of our method, and are omitted simply for the sake of brevity; for a flavor
of these results, the reader is referred to \cite{JLB:TezaurWCCM}.

In future work, we plan to explore the following research directions: 
\begin{itemize}
\item \textit{Verifying that similar results can be obtained via explicit-explicit and implicit-explicit Schwarz couplings, and in the case when disparate time-steps are used to advance 
the solution in different subdomains.}  Preliminary testing not reported here suggests that the same conclusions hold in the case of couplings involving disparate time-integrators/time-steps.
    \item \textit{An extension of the proposed coupling framework to 2D and 3D problems.}  A multi-dimensional implementation will require the development of transfer operators for defining the Schwarz transmission boundary conditions.  We plan to investigate
the usage of the Compadre toolkit \cite{JLB:Compadre} for this task.  It is expected that the benefits of hyper-reduction will be far greater in 2D and 3D than in 1D.  
    \item \textit{Developing error indicator-based approaches for determining ``optimal" domain decompositions and ROM/FOM placement.}  We plan to build on the work of Bergmann \textit{et al.} \cite{JLB:Bergmann2018}.
    \item \textit{Implementing and testing the additive Schwarz variant of the proposed coupling method.}  Please see Remark 3 and \cite{JLB:Gander2008} for more information on this variant of the method.  Preliminary studies suggest that additive Schwarz-based coupling has the potential to yield  coupled models which can achieve speed-ups over a single-domain FOM when parallelized over the number of subdomains.
    \item \textit{Examining snapshot collection strategies that do not require simulating the coupled
high-fidelity model over the the entire domain.}  Ideas such as oversampling \cite{JLB:Smetana2022} will be explored towards the goal of designing a workflow that can reuse existing modular codes for individual physics in an effort to perform  minimally-intrusive modular couplings.
    \item \textit{Extending the proposed approach to  coupling scenarios that enable ``on-the-fly" FOM-ROM switching and ROM adaptation.}  A nice starting point for this research direction is the work of Corigliano \textit{et al.}  \cite{JLB:Corigliano1}.
	\item \textit{Performing an analysis of the method's theoretical convergence properties.}  We plan to leverage some of our past work, which analyzed the method's convergence for FOM-FOM coupling in solid mechanics \cite{JLB:mota2017schwarz, JLB:mota2022schwarz}.
	\item \textit{Applying the proposed coupling methodology to problems other problems, including multi-physics problems.}  We are particularly interested in 
	problems involving FSI.  We have begun to move in this direction by starting to explore multi-material coupling using Schwarz, as a proxy for a multi-physics problem,  
	\item \textit{Extending the proposed framework to enable the seamless coupling of other types of data-driven models, e.g., PINNs.}  PINNs are well-suited for Schwarz-based couplings, as they have the notion of boundary conditions, which can be incorporated through the loss function being minimized \cite{JLB:LiD3M,JLB:LiDeepDDM}.
\end{itemize}

\bibliographystyle{siam}
\bibliography{JoshuaBarnett.bib}


\end{document}